\documentclass[twoside,german,10pt,a4]{amsart}
\usepackage{amssymb,amsmath}

\newtheorem{theo}{\sc Theorem}

 \newtheorem{prop}{\sc Proposition}

\def\det{{\rm det\,}}

\def\Res{\mathop{\rm res}}
\def\N{\mathbb{N}}

\def\C{\mathbb{C}}
\def\Z{\mathbb{Z}}

\def\Re{{\rm Re} \,}
\def\Im{{\rm Im} \,}

\def\2int{\mathop{\int\!\!\!\int}}

\def\kzu1#1{\buildrel #1 \over \longrightarrow}
\def\mod{{\rm mod\,}}

\def\dist{{\rm dist\,}}

\def\ds{\displaystyle}
\def\ts{\textstyle}

\def\und{\quad{\rm and}\quad}

\def\Y{\mathfrak{Y}}
\def\barY{\overline\Y}\def\barY{\widetilde\Y}

\def\u{\mathfrak{u}}

\def\y{\mathfrak{y}}

\def\z{\mathfrak{z}}

\def\c{\mathfrak{c}}

\def\f{\mathbf{f}}
\def\f{\mathfrak{f}}

\def\h{\mathbf{h}}

\def\n{\mathfrak{n}}
\def\p{\mathfrak{p}}
\def\q{\mathfrak{q}}

\def\0{\mathbf{0}}
\def\1{\mathbf{1}}

\def\P{\mathfrak{P}}

\def\CL{\mathfrak{C}}
\def\be#1{\begin{equation}\label{#1}}
\def\ee{\end{equation}}
\def\ende{{\blacksquare}}
\def\Ende{$\ende$}
\def\schema#1#2#3#4{\begin{array}{c|c}#2&#1\cr\hline #3&#4\end{array}}

\def\horizontal#1{\begin{array}{cccc}#1\!\!\!&#1\!\!\!&#1\!\!\!&#1\end{array}}
\def\vertikal#1{\begin{array}{c}#1\cr#1\cr#1\cr#1\end{array}}
\def\symbol{\diamond}

\begin{document}
\title{\sc An old new class of meromorphic functions}
\author{\sc Norbert Steinmetz}
\maketitle

\bigskip

{\small\begin{center} {\sc Abstract}\end{center}
\begin{quote} Based on the so-called re-scaling method, we will give a detailed description of the solutions to the Hamiltonian system (\ref{Hsystem}) below, which was discovered only recently by Kecker~\cite{TK1,TK2}, and is strongly related to Painlev\'e's fourth differential equation. In particular, the problem to determine those fourth Painlev\'e transcendents with positive Nevanlinna deficiency $\delta(0,w)$, is completely resolved.
\end{quote}


}

\section{Introduction}

The Hamiltonian system
\be{Hsystem} p'=-q^2-zp-\alpha,~q'=p^2+zq+\beta\ee
with time-dependent Hamiltonian
\be{hamiltonian} H(z,p,q)=\ts\frac13(p^3+q^3)+zpq+\beta p+\alpha q\ee
has been discovered by Kecker~\cite{TK1,TK2} when questing for systems having the so-called (pseudo) {\it Painlev\'e property}. Roughly speaking, this means that the system has ``no movable singularities but poles'', or, more precisely, that every solution admits unrestricted analytic continuation in the plane except for fixed singularities
(which do not occur here); and ``pseudo'' means that, in addition, also moving algebraic singularities are admitted. It is not hard to verify, although more elaborate to discover, that $w=p+q-z$ solves the differential equation
\be{PIVmod}2ww''=w'^2-w^4-4zw^3-(2\alpha+2\beta+3z^2)w^2-(\alpha-\beta+1)^2\,.\ee
Kecker~\cite{TK1} found some implicit second order equation for $q$.
Obviously, equation (\ref{PIVmod}) is closely related to Painlev\'e's fourth equation
\be{PIV}2yy''=y'^2+3y^4+8zy^3+4(z^2-\hat\alpha)y^2+2\hat\beta\ee
with parameters $\hat\alpha=\frac{i}{\sqrt 3}(\alpha+\beta)$ and $\hat\beta=\frac29(\alpha-\beta+1)^2;$ one has just to consider $y(z)=aw(bz)$ with $b=\sqrt[4]{-\frac43}$ and $a=-\frac12b^3$.
Since also $\tilde w=\omega p+\bar\omega q-z$ ($\omega^3=1$) is a solution to (\ref{PIVmod}) with $(\alpha,\beta)$ replaced by $(\omega\alpha,\bar\omega\beta)$,
$p$ and $q$ may be re-discovered:
\be{RETOUR}(\omega-\bar\omega)p=\tilde w-\bar\omega w-(\bar\omega-1)z\und (\bar\omega-\omega)q=\tilde w-\omega w-(\omega-1)z.\ee
This could be the end of the story. The ``old new'' transcendents $p$ and $q$, however, have so many interesting properties that it seems justified to study them in their own right.\medskip

\section{Notation and simple properties}
By $\Lambda$ we denote the set of non-zero poles of $p$. It is easily seen that the poles $\lambda$ are simple, and $p$ and $q$ have residues $\varrho$ and $-\bar\varrho$, respectively,
where $\varrho$ is some third root of unity. Thus $\Lambda$ is  divided in a natural way into three disjoint sets  $\Lambda^\varrho$, $\varrho^3=1$.
With $\Lambda$ we associate the sets
\be{TRIANGLE}\triangle_\delta(\lambda)=\{z:|z-\lambda|<\delta|\lambda|^{-1}\}\und \Lambda_\delta=\bigcup_{\lambda\in\Lambda}\triangle_\delta(\lambda).\ee
If it is clear which solution  $(p,q)$ is under consideration, we will also write $H(z)=H(z,p(z),q(z)).$
The following holds at poles in $\Lambda^1$:
\be{POLES}\begin{array}{rcl}
p(z)&=&\phantom{-}(z-\lambda)^{-1}+\ts\frac12\lambda+\big(1+\frac13(\alpha-2\beta)-\frac14{\lambda^2}\big)(z-\lambda)\cr
&&\quad +\big(\h-\big(\frac58+\frac14(\alpha-\beta)\big)\lambda\big)(z-\lambda)^2+\cdots\cr
q(z)&=&-(z-\lambda)^{-1}+\ts\frac12\lambda+\big(1+\frac13(2\alpha-\beta)+\frac14{\lambda^2}\big)(z-\lambda)\cr
&&\quad +\big(\h+\big(\frac58+\frac14(\alpha-\beta)\big)\lambda\big)(z-\lambda)^2+\cdots\end{array}\ee
\be{POLESH} H(z)=\ts (z-\lambda)^{-1}+[2\h +\frac 13\lambda^3+\frac12(\alpha+\beta)\lambda]+
\big[\frac13(\alpha+\beta)+\frac34\lambda^2\big](z-\lambda)+\cdots\ee
\be{POLESDiv}p(z)+q(z)-z=(1+\alpha-\beta)(z-\lambda)+2\h(z-\lambda)^2+\cdots\ee
The coefficient $\h=\h(\lambda)$ remains undetermined and {\it free}. Actually $\lambda$ and $\h$ may be prescribed to determine a unique solution.
To determine the Laurent series at poles with residue $\varrho$ ($\varrho^3=1$), we replace $p$ and $q$ by $x=\bar\varrho p$ and $y=\varrho q$, respectively. Then $x$ and $y$ satisfy $x'=-y^2-zx-\bar\varrho\alpha,$ $ y'=x^2+zy+\varrho\beta,$
and we thus obtain the Laurent developments for $p=\varrho x$ and $q=\bar\varrho y$ by replacing
$(p,q,\alpha,\beta)$ in (\ref{POLES}) by $(x,y,\bar\varrho\alpha,\varrho\beta)$:
\be{pqomega}\begin{array}{rcl}
p(z)&=&\phantom{-}\varrho(z-\lambda)^{-1}+\frac12\varrho \lambda+\big(\varrho+\frac13(\alpha-2\bar\varrho\beta)-\frac14{\varrho}\lambda^2\big)(z-\lambda)+\cdots\cr
q(z)&=&-\bar\varrho(z-\lambda)^{-1}+\frac12{\bar\varrho}\lambda+\big(\bar\varrho+\frac13(2\varrho\alpha-\beta)+\frac14{\bar\varrho}\lambda^2\big)(z-\lambda)
+\cdots\end{array}\ee
In particular, this means that the development (\ref{POLESH})
remains valid at poles with residue $\varrho$, provided $\alpha$ and $\beta$ are  replaced by $\bar\varrho\alpha$ and $\varrho\beta$, respectively.
Hence $g(z)=\exp({\ts \int} H(z)\,dz)$
is an entire function that has simple zeros at the poles of $p$, and no others.

\subsection{B\"acklund transformations} Trivial B\"acklund transformations are $(\omega^3=1$ arbitrary)
\be{TBT}\begin{array}{rl}
{\sf M}_\omega:&\left\{\begin{array}{c}\tilde p(z)=\bar\omega p(z),~\tilde q(z)=\omega q(z)\cr
(\alpha,\beta)\mapsto(\bar\omega\alpha,\omega\beta)\end{array}\right.\cr
{\sf R}:& \left\{\begin{array}{c}\tilde p(z)=-iq(iz),~\tilde q(z)=-ip(iz)\cr
(\alpha,\beta)\mapsto(-\beta,-\alpha)\end{array}\right.\cr
\overline{\sf R}:& \left\{\begin{array}{c}\tilde p(z)=\overline{p(\bar z)},~\tilde q(z)=\overline{q(\bar z)}\cr
(\alpha,\beta)\mapsto(\bar\alpha,\bar\beta)\end{array}\right.\end{array}\ee
Nontrivial B\"acklund transformations were found by Kecker~\cite{TK1}
\be{NTB}{\sf B}_\omega:~
\left\{\begin{array}{c}
\tilde p(z)=\ds p(z)-\bar\omega\frac{\omega\alpha-\bar\omega\beta+1}{\omega p(z)+\bar\omega q(z)-z}\cr
\tilde q(z)=\ds q(z)+\omega\frac{\omega\alpha-\bar\omega\beta+1}{\omega p(z)+\bar\omega q(z)-z}\cr
(\alpha,\beta)\mapsto(\omega\beta-\bar\omega,\bar\omega\alpha+\omega)\end{array}\right.\ee
(as long as the denominator does not vanish identically). By {\it B\"acklund transformation} we mean any repeated
application of the above transformations. B\"acklund transformations act on pairs $(p,q)$, components $p$ and $q$, and also parameters $(\alpha,\beta)$.

Under the B\"acklund transformation ${\sf B}_\omega$ the residues  $\varrho=\Res_\lambda p$ change as follows:
\be{CHANGEofRESIDUES}
\begin{array}{rl}
\varrho\ne\bar\omega:&\Res_\lambda{\sf B}_\omega p=\varrho\cr
\varrho=\bar\omega:& {\sf B}_\omega p {\rm ~is~regular~at~}\lambda {\rm ~and~}\cr
&\omega p+\bar\omega q-z|_{z=\lambda}=0\cr
\left.\begin{array}{r}\omega p(\tilde\lambda)+\bar\omega q(\tilde\lambda)-\tilde\lambda=0\cr
(\tilde\lambda {\rm~not~a~pole~of~} p)\cr
\omega p'(\tilde\lambda)+\bar\omega q'(\tilde\lambda)=
-\omega\alpha+\bar\omega\beta\end{array}\right\}
:& \Res_{\tilde\lambda}{\sf B}_\omega p=\bar\omega\end{array}\ee

\section{Rescaling}

\subsection{Yosida functions}\index{Yosidafunctions}
Let $a$ and $b>-1$ be real parameters. By definition, the class $\barY_{a,b}$ consists of all meromorphic functions $f$ such that
the family $(f_\kappa)_{|\kappa|>1}$ of functions
\be{RESCALE}f_\kappa(\z)=\kappa^{-a}f(\kappa+\kappa^{-b}\z)\ee
is normal on $\C$ in the sense of Montel, and all limit functions  $\mathfrak{f}=\lim\limits_{\kappa_n\to\infty}f_{\kappa_n}$ are $\not\equiv\infty$, at least one of them being non-constant.
If, in addition, {\it all} limit functions are non-constant, then $f$ is said to belong to the {\it Yosida class}\index{Yosida class}
$\Y_{a,b}$. The functions of class $\Y_{0,0}$ were introduced by Yosida~\cite{Yosida2}, and for arbitrary real parameters by the author~\cite{NStYosida}.
The class $\Y_{0,0}$ is universal in the sense that it contains all limit functions $\f=\lim\limits_{\kappa_n\to\infty}f_{\kappa_n}$
for $f\in\Y_{a,b}$. The functions $f\in\barY_{a,b}$  have
striking properties, for example they satisfy $f^\sharp(z)=O(|z|^{|a|+b})$, $T(r,f)=O(r^{2+2b})$, and $m(r,f)=O(\log r)$,
and even $T(r,f)\asymp r^{2+2b}$ and $m(r,1/f')=O(\log r)$ if $f\in\Y_{a,b}$. For notations and results in Nevanlinna Theory the reader is referred to Hayman's monograph~\cite{Hayman}.

\subsection{An application to the system (\ref{Hsystem})} We quote from Shimomura's paper \cite{Shimomura1}, section 5.2, and also from section 6 in the author's paper \cite{NStPainleveOrder} the following facts about
the solutions to (\ref{PIV}) and, obviously, also to the functions $w=p+q-z$ solving (\ref{PIVmod}) and to $w=\bar\varrho p+\varrho q-z$ ($\varrho^3=1$)
which solve similar equations:
\begin{itemize}
\item For $\delta>0$ sufficiently small, the discs $\triangle_\delta(\lambda)$ ($\lambda\in\Lambda$) are mutually disjoint;
\item $w=O(|z|)$ as $z\to\infty$ outside $\Lambda_\delta$ (for the definition of $\triangle_\delta$ and $\Lambda_\delta$ cf. (\ref{TRIANGLE})).
\end{itemize}

It follows from (\ref{RETOUR}) that the second condition is equivalent to
\be{pqest}|p|+|q|=O(|z|)\quad(z\to\infty, z\notin\Lambda_\delta),\ee
while (\ref{Hsystem}) implies
\be{pqprimeest}|p'|+|q'|=O(|z|^2)\quad(z\to\infty, z\notin\Lambda_\delta).\ee
If $\lambda$ is any pole with residue $\Res_\lambda p=\varrho$, then it follows from (\ref{POLES}) that
$f=\varrho p'+p^2-\varrho zp$ is regular on $\triangle_\delta(\lambda)$, and $f(z)=O(|z|^2)$ on $\partial\triangle_\delta(\lambda)$ continues to hold on $\triangle_\delta(\lambda)$
by the Maximum Principle. This implies $|p'|=O(|z|^2+|p|^2)$ on $\triangle_\delta(\lambda)$, and on combination with (\ref{pqprimeest}) this gives
\be{pqNormal}\frac{|p'|}{|z|^2+|p|^2}+\frac{|q'|}{|z|^2+|q|^2}=O(1)\quad(z\to\infty)~{\rm without~restriction}.\ee

From (\ref{pqNormal}) easily follows:

\begin{theo}\label{pqnormaluty}The solutions $p$ and $q$ to equation (\ref{Hsystem}) belong to the class $\barY_{1,1}$.
The limit functions $\p=\lim\limits_{\kappa_n\to\infty}p_{\kappa_n}$ and $\q=\lim\limits_{\kappa_n\to\infty}q_{\kappa_n}$
solve
\be{LIMIT}\p'=-\q^2-\p,~\q'=\p^2+\q,~\ts\frac13(\p^3+\q^3)+\p\q=\c\ee
for some constant $\c$ depending on $(\kappa_n)$.\end{theo}

\proof  Set  $f(z)=p(z)/z$, $z=\kappa+\kappa^{-1}\z$ and $f_\kappa(\z)=f(z)$. Then (\ref{pqNormal}) implies
$$f^\sharp(z)\le\frac{|p'(z)||z|^{-1}}{1+|p(z)|^2|z|^{-2}}+\frac{|p(z)||z|^{-2}}{1+|p(z)|^2|z|^{-2}}=|z|\frac{|p'(z)|}{|z|^2+|p(z)|^2}+O(1)=O(1+|z|)$$
and
$$f^\sharp_\kappa(\z)=|\kappa|^{-1}f^\sharp(\kappa+\kappa^{-1}\z)=O(1+|\kappa|^{-2}|\z|)=O(|\z|)\quad(\z\to\infty).$$
By Marty's well-known criterion (see Ahlfors~\cite{Ahlfors1}), the family $(f_\kappa)_{|\kappa|>1}$ is normal on $\C$, and so is the family of functions $p_\kappa(\z)=(1+\kappa^{-2}\z)f_\kappa(\z).$
The same is true for the family $(q_\kappa)$. The limit functions
are finite since $\p$ and $\q$ have simple poles at $\z=0$
with residues $\varrho$ and $-\bar\varrho$, respectively, if  $\lim\limits_{\kappa_n\to\infty}|\kappa_n|\dist(\kappa_n,\Lambda)=0$, while $\p(0)$ and $\q(0)$ are finite if $\liminf\limits_{\kappa_n\to\infty}|\kappa_n|\dist(\kappa_n,\Lambda)>0$. Obviously, $\p$ and $\q$ solve (\ref{LIMIT}) with
\be{LimitcNonpole}\c=\lim\limits_{n\to\infty} \kappa_n^{-3} H(\kappa_n)\quad{\rm if}~\liminf\limits_{\kappa_n\to\infty}|\kappa_n|\dist(\kappa_n,\Lambda)>0\ee
and
\be{Limitcpole}\c=\lim\limits_{n\to\infty} \lambda_n^{-3}(2\h(\lambda_n)+{\ts\frac 13}\lambda_n^3)\quad{\rm if}~\lim\limits_{\kappa_n\to\infty}|\kappa_n||\kappa_n-\lambda_n|=0\ee
for some sequence of poles $(\lambda_n)$. \Ende

\medskip{\sc Remark.} Constant limit functions are $(0,0)$ and $(-\omega,-\bar\omega)$ with $\omega^3=1$. They correspond to $\c=0$ and $\c=\frac13$, respectively.

\subsection{The cluster set}
By definition, the {\it cluster set} $\CL(p,q)\subset\C$ of any non-trivial solution $(p,q)$ to (\ref{Hsystem}) consists of all limits (\ref{LimitcNonpole}).

\begin{theo}\label{CLSET}The cluster set $\CL(p,q)$ is closed, bounded, and connected, and contains all limits
(\ref{Limitcpole}).\end{theo}

\proof For $\delta>0$ sufficiently small, the closed discs $\bar\triangle_\delta(\lambda)$ about the poles $\lambda\ne 0$ are mutually disjoint, hence any two points $a,b\in D_\delta=\C\setminus\bigcup_{\lambda\in\Lambda}\bar\triangle_\delta(\lambda)$  may be joined by a curve that is contained in $D_\delta\cap\{z:|z|\ge \min\{|a|,|b|\}\}$. We denote the corresponding cluster set of $z^{-3} H(z)$ as $z\to\infty$ on $D_\delta$ by $\CL_\delta(p,q)$, and
note that $\CL(p,q)=\bigcup_{\delta>0}\CL_\delta(p,q)$. The cluster sets $\CL_\delta(p,q)$ are closed, bounded
in $\C$ since $z^{-3} H(z)$ is uniformly bounded on $D_\delta$, and connected by the special property of $D_\delta$.
Since $\CL_\eta(p,q)\supset\CL_\delta(p,q)$ for $0<\eta<\delta,$ it remains to show that
\be{CetasubCdelta}\CL_\eta(p,q)\subset\CL_\delta(p,q)\quad(0<\eta<\delta {\rm ~sufficiently~small}),\ee
and that the limits (\ref{Limitcpole}) belong to $\CL_\delta(p,q)$. If
$$\eta\le |\kappa_n|\dist(\kappa_n,\Lambda)=|\kappa_n||\kappa_n-\lambda_n|<\delta$$
holds for some $\lambda_n\in\Lambda$, and if $\p=\lim\limits_{\kappa_n\to\infty}p_{\kappa_n}$, $\q=\lim\limits_{\kappa_n\to\infty}q_{\kappa_n}$ and $\c=\lim\limits_{\kappa_n\to\infty}\kappa_n^{-3} H(\kappa_n)$ exist, we replace $\kappa_n$ by $\lambda_n$
with the following effect: from $\lambda_n=\kappa_n+\kappa_n^{-1}\z_n$ with $(\z_n)$ bounded, hence $\z_n\to\z_0$ as we may assume, it follows that
$\lambda_n^{-1}=(1+o(1))\kappa_n^{-1}$,
$p_{\kappa_n}(\z)=(1+o(1))p_{\lambda_n}(\z-\z_0+o(1))$ and $ q_{\kappa_n}(\z)=(1+o(1))q_{\lambda_n}(\z-\z_0+o(1)),$
hence
$$\hat\p(\z)=\lim_{\lambda_n\to\infty}p_{\lambda_n}(\z)=\p(\z+\z_0),~ \hat\q(\z)=\lim_{\lambda_n\to\infty}q_{\lambda_n}(\z)=\q(\z+\z_0)$$
and
$${\ts\frac 13}(\hat\p^3+\hat\q^3)+\hat\p\hat\q=\c=\lim\limits_{n\to\infty} \lambda_n^{-3}(2\h(\lambda_n)+{\ts\frac 13}\lambda_n^3).$$

Finally, if we start with sequences $(p_{\lambda_n})$ and $(q_{\lambda_n})$ with $\lambda_n\in\Lambda$, we may as well consider $(p_{\kappa_n})$ and $(q_{\kappa_n})$
with $|\kappa_n||\kappa_n-\lambda_n|=\delta$, hence
$$\c=\lim\limits_{\lambda_n\to\infty}\lambda_n^{-3}(2\h(\lambda_n)+{\ts\frac 13}\lambda_n^{3})=\lim_{\kappa_n\to\infty}\kappa_n^{-3} H(\kappa_n)\in\CL_\delta(p,q).$$
Combining both arguments this also proves (\ref{CetasubCdelta}). \Ende

\subsection{An algebraic curve}
The {\it algebraic curve}
\be{algcurve}\ts\frac13(u^3+ v^3)+uv=c\ee
\begin{itemize}\item
is {\it reducible} if $c=\frac 13$: $u^3+v^3+3uv-1=\prod\limits_{\varrho^3=1}(v+\varrho u-\bar\varrho)$,
\item has {\it genus} zero if $c=0$, and
\item has {\it genus} one otherwise.\end{itemize}
In the first case $(c=\frac13)$ the corresponding Hamiltonian system
\be{hamsys}u'=-v^2-u,~v'=u^2+v\ee
has solutions given by $\ds\frac{u+\omega}{u+\bar\omega}=e^{i\sqrt 3 t}$ ($\omega=\frac12(-1+i\sqrt 3)$) and $v=1-u$, and similar expressions if $v=\bar\varrho-\varrho u$ and $v=\varrho-\bar\varrho u$, respectively. It is just important to know that the poles form a ${2\pi}/{\sqrt 3}-$periodic sequence with fixed residue.
If the genus is zero $(c=0)$, the system (\ref{hamsys}) is solved by $u=\ds\frac{-3\omega e^{2t}}{e^{3t}+1}$ and $v=\ds\frac{-3\bar\omega e^{t}}{e^{3t}+1}$,
with poles forming a $2\pi i/3-$periodic sequence, this time with alternating residues.
In the remaining cases, the solutions to (\ref{hamsys}) are elliptic functions parametrising the curve (\ref{algcurve}); they have elliptic order three with corresponding lattice $L_c$ depending only on $c$.

\section{Value distribution}\label{VD}

\subsection{The order of growth} We note as a corollary to Theorem~\ref{pqnormaluty} the following estimate:

\begin{theo}\label{ORDER}The solutions to equation (\ref{Hsystem}) have order of growth at most four:
\be{Trpq}T(r,p)+T(r,q)=O(r^4).\ee\end{theo}

\proof From $f^\sharp(z)=O(|z|)$ and $T(r,p)=T(r,f)+O(\log r)$ for $f(z)=p(z)/z$ we easily obtain (\ref{Trpq}) when using the Ahlfors-Shimizu form of the Nevanlinna characteristic,
$$T(r,f)=\int_0^rA(t)\frac{dt}t\quad{\rm with}\quad A(t)=\frac1\pi\int_{|z|<t}f^\sharp(z)^2\,dx\,dy.\quad\ende$$

\medskip {\sc Remark.} This result corresponds to the meanwhile well-established estimate \cite{Shimomura1,NStPainleveOrder} for the order of growth of the fourth Painlev\'e transcendents.

\subsection{Solutions of maximal growth} Suppose that the cluster set $\CL(p,q)$ contains some parameter $\lim\limits_{ \kappa_n\to\infty} H( \kappa_n) \kappa_n^{-3}\ne 0,\frac13$.
Since $H(z)z^{-3}$ is  varying slowly, this following from
$$\frac d{dz} H(z)z^{-3}=-3 H(z)z^{-4}+p(z)q(z)z^{-3}=O(|z|^{-1})\quad{\rm outside~} \Lambda_\delta,$$
given $\epsilon>0$ sufficiently small there exists $\eta>0$ such that the re-scaling procedure for $|\tilde \kappa_n- \kappa_n|<\eta| \kappa_n|$ leads to limit functions
$\p$ and $\q$ with
$$\p^3+\q^3+3\p\q=3\c\ts\quad(\min\{|\c|,|\c-\frac13|\}\ge\epsilon,~|\c|\le 1/\epsilon).$$
Now $\p$ and $\q$ are elliptic functions with fundamental parallelogram $P_\c$, whose diameter and area is uniformly bounded and bounded away from zero.
It is thus easily deduced that the disc $|z- \kappa_n|<\eta| \kappa_n|$ contains at least const$| \kappa_n|^4$ poles of $p$ (it is almost the same to say that the disc $|z|<r$ contains $\sim \pi r^2$ lattice points $m+in$), hence
$$n(2| \kappa_n|,p)\ge const|\kappa_n|^4$$
holds. On combination with $n(r,\Lambda)=O(r^4)$ we thus have:

\begin{theo}\label{ORDERDICHO}Let $(p,q)$ be any solution to equation (\ref{Hsystem}). Then $\CL(p,q)\not\subset\{0,{\ts\frac13}\}$
implies $n(r_k,\Lambda)\asymp r_k^4$, at least on some sequence $r_k\to\infty$.\end{theo}

\subsection{The distribution of residues}Leaving the rational solutions aside we henceforth will consider only transcendental solutions.
From
$$q^2=-p'-zp-\alpha=-p\,(z+p'/p)-\alpha$$
it follows by the usual rules of Nevanlinna Theory that $2m(r,q)\le m(r,p)+O(\log r)$, and in the same manner
$2m(r,p)\le m(r,q)+O(\log r)$ is obtained, hence $m(r,p)+m(r,q)=O(\log r)$,
$T(r,p)=N(r,\Lambda)+O(\log r),$ and $T(r,q)=N(r,\Lambda)+O(\log r)$
hold; here $N(r,\Lambda)$ denotes the common Nevanlinna counting function of poles.
Cauchy's Residue Theorem yields
\be{INT}\frac1{2\pi i}\int_{C_r}p(z)\,dz=n(r,\Lambda^1)+\varrho n(r,\Lambda^\varrho)+\bar\varrho n(r,\Lambda^{\bar\varrho})\quad(\varrho=\ts\frac12(-1+i\sqrt 3)),\ee
provided the circle $C_r:|z|=r$ intersects no pole. If $C_r$ intersects the disc $\triangle_\delta(\lambda)$, we replace the arc
$C_r\cap\triangle_\delta(\lambda)$ by the sub-arc of $\partial\triangle_\delta(\lambda)$ outside $C_r$ if $|\lambda|\le r$, and inside $C_r$ otherwise.
This way we obtain a simple closed curve $\Gamma_r$ such that
$$|p(z)|+|q(z)|=O(|z|)\quad{\rm holds~on~}\Gamma_r,$$
without changing the integral (\ref{INT}). The length of the part of $\Gamma_r$ restricted to any sector of central angle $\Theta$
is $O(r\Theta)$, hence the integral is $O(r^2)$. Taking real and imaginary parts, we obtain
\be{REALIMAGPART}\begin{array}{rcl}
n(r,\Lambda^1)-\frac12(n(r,\Lambda^\varrho)+n(r,\Lambda^{\bar\varrho}))&=&O(r^2)\phantom{.}\cr
\frac{\sqrt 3}2(n(r,\Lambda^{\varrho})-n(r,\Lambda^{\bar\varrho}))&=&O(r^2)\end{array}\quad(\varrho={\ts\frac12(-1+i\sqrt 3)}).\ee

\begin{theo}\label{EQUALN}For any transcendental solution $(p,q)$ to equation (\ref{Hsystem}), the counting functions $n(r,\Lambda^\varrho)$
are equal up to a term $O(r^2)$, that is, we have
$$n(r,\Lambda^\varrho)={\ts\frac 13}n(r,\Lambda)+O(r^2)\quad(\varrho^3=1).$$
In particular, if one kind of poles is missing, $p$ and $q$ have order of growth at most two:
$T(r,p)+T(r,q)=O(r^2).$\end{theo}

\subsection{Strings of poles}Re-scaling along any sequence of poles with corresponding limit $\c\in\{0,{\ts\frac13}\}$ leads to the following situation: given $\epsilon>0$ and $R>0$ there exists $r_0>0$,
such that for any pole $\lambda_0$ in $|z|>r_0$, the disc $\triangle_R(\lambda)$ contains the poles
$$\lambda_k=\lambda_0+k(\varpi+\epsilon_k)\lambda_0^{-1}\quad(|\epsilon_k|<\epsilon,~-k_1\le k\le k_2)$$
with $\varpi=2\pi i/3$ and $\varpi=2\pi/\sqrt 3$, respectively, and no others
($k_1=k_2=[R/|\varpi|]$ if $R$ is large and not an integer multiple of $|\varpi|$).
In other words, for $\CL(p,q)\subset\{0,\frac13\}$, each pole $\lambda$ with $|\lambda|$ sufficiently large belongs to some unique sequence $(\lambda_n)$, called {\it string of poles}; it satisfies the approximative recursion
\be{lambdarek}\lambda_{n+1}=\lambda_n\pm(\varpi+o(1))\lambda_n^{-1}.\ee
Setting $\sigma_n=\lambda_n^2$ we obtain $\sigma_{n+1}=\sigma_n\pm 2\varpi+o(1)$, $\sigma_n=2n\varpi+o(n)$,
$$\lambda_n=\sqrt{2n|\varpi|}(1+o(1))\und \arg\lambda_n=\ts\frac12\arg\varpi+o(1)~\mod\frac\pi 2;$$
the counting function of $\{\lambda_n\}$ is given by $n(r,\{\lambda_n\})\sim\frac{r^2}{2|\varpi|}.$
We note that the estimate $n(r,\Lambda)=O(r^4)$ now can be completed:

\begin{theo}\label{UNTEN}Any transcendental solution $(p,q)$ to (\ref{Hsystem}) satisfies
$$r^2=O(n(r,\Lambda))\und n(r,\Lambda)=O(r^4).$$\end{theo}

\proof This is obvious if $\Lambda$ contains a full string. In any case the re-scaling method shows that to any pole $\lambda$ of sufficiently large modulus there exists some pole $\tilde\lambda$
satisfying $|\lambda|<|\tilde\lambda|<|\lambda|+O(|\lambda|^{-1})$, this implying $r^2=O(n(r,\Lambda))$. \Ende

\section{Asymptotic expansions and pole-free sectors}

\subsection{Pole-free sectors}Let $(p,q)$ be any transcendental solution to (\ref{Hsystem}). If $p$ has no poles on some sector $S:|\arg z-\hat\theta|<\theta,~|z|>r_0$, the re-scaling procedure with $\kappa_n\in S_\delta:|\arg z-\hat\theta|<\theta-\delta$
leads to solutions $(\p,\q)$ without poles, hence to constant solutions. The possible constants are $(0,0)$ and $(-\tau,-\bar\tau)$, again with
$\tau^3=1$, hence we have either $p(z)=o(|z|)$ and $q(z)=o(|z|)$, or else $p(z)=-\tau z+o(|z|)$ and $q(z)=-\bar\tau z+o(|z|)$ as $z\to\infty$ on $S$.
In each case we will prove that this leads to asymptotic expansions on certain pole-free sectors.

\subsection{Asymptotics on pole-free sectors with $\mathbf{\CL(p,q)=\{0\}}$}In the first case the following holds in more generality:

\begin{theo}\label{ASY0}Suppose that $p(z)=o(|z|)$ and $q(z)=o(|z|)$ hold as $z\to\infty$ on some ray $\hat\sigma:\arg z=\hat\theta\not\equiv {\ts\frac\pi 4}~\mod {\ts\frac\pi 2}$. Then $p$, $q$, and $H$ have asymptotic expansions
\be{ASYMP0}\begin{array}{rcl}
p(z)&\sim&\ds-\frac\alpha z-\frac{\alpha+\beta^2}{z^3}-\frac{3\alpha+\beta^2+2\alpha^2\beta}{z^5}+\cdots\cr
q(z)&\sim&\ds-\frac\beta z+\frac{\beta-\alpha^2}{z^3}+\frac{3\beta-\alpha^2-2\alpha\beta^2}{z^5}+\cdots\cr
 H(z)&\sim&\ds-\frac{\alpha\beta}z-\frac{\alpha^3+\beta^3}{3z^3}+\cdots\qquad(z\to\infty)
\end{array}\ee
on the sector $\Sigma_\nu:|\arg z-\nu{\ts\frac\pi 2}|<{\ts\frac\pi 4}$ that contains the ray $\hat\sigma$.\end{theo}

\proof From our hypothesis it also follows that $H(z)=o(|z|^3)$ holds along $\hat\sigma$. Re-scaling along any sequence $(\kappa_n)$ on  $\hat\sigma$ yields limit functions $\p$ and $\q$
satisfying
$$\p'=-\q^2-\p,~\q'=\p^2+\q,~\p(0)=\q(0)=0, ~{\rm and}~ \p^3+\q^3+3\p\q=0,$$
hence $\p=\q\equiv 0$.
Thus there are pole-free discs $|z-re^{i\hat\theta}|<r^{-1}\kappa(r)$ such that $\kappa(r)\to\infty$ as $r\to\infty$.
To identify the maximal pole-free sector that contains $\hat\sigma$ (if any), we start with $r_0>0$ sufficiently large and define the sequence $(r_n)$ inductively by
$r_{n+1}=r_n+4 r_n^{-1}.$ By $\theta_n$ we denote the largest number such that
$$A_n=\{z:r_n\le|z|\le r_{n+1},~\hat\theta\le\arg z<\theta_n\}$$
contains no pole of $p$, noting that $r_n\theta_n\sim \kappa(r_n)\to\infty$. We may assume that there exists some pole $\tilde\lambda_0$ on $\partial A_0$, hence $\arg\tilde\lambda_0=\theta_0$.
The same is true at least for some
sub-sequence $\partial A_{n_k}$. Re-scaling along the sequence $(\lambda_{n_k})$ then yields limit functions $\tilde\p$ and $\tilde\q$ that have simple poles at $\z=0$
and satisfy $\tilde\p^3+\tilde\q^3+3\tilde\p\tilde\q=3\c$. Since $A_{n_k}$ is large with respect to the metric $ds=|z||dz|$ and contains no pole, $\tilde\p$ and $\tilde\q$ cannot be elliptic functions, and this
and $\lim\limits_{z\to\infty} H(z)z^{-3}=0$ on $\bigcup_n A_n$ implies $\c=0$, hence the algebraic curve (\ref{algcurve}) has genus zero. It follows from Hurwitz' Theorem
on the poles of limit functions that to each pole $\tilde\lambda_{n_k}$ there exist five poles
$$z_{k,\nu}=\tilde\lambda_{n_k}+\nu(\varpi+o(1))\tilde\lambda_{n_k}^{-1}\quad(-2\le \nu\le 2,~\varpi=\pm 2\pi i/3)$$
on $|z-\tilde\lambda_n|<5|\tilde\lambda_n|^{-1}$, and no others. Since $z_{k,\pm 2}$ does not belong to the annulus $r_{n_k}\le|z|\le r_{n_k+1}$,
it follows that $A_{n_k-1}$ and $A_{n_k+1}$, hence {\it each} $A_n$ contains some pole $\tilde\lambda_n$ on its boundary, and $(\tilde\lambda_n)$
is a sub-sequence of some sequence $(\lambda_n)$ of poles satisfying the approximate recursion
(\ref{lambdarek}). Thus $\arg\lambda_n$ and $\theta_n$ approach  $(2\nu+1){\ts\frac\pi 4}>\hat\theta$ for some $\nu$. The same argument applies to the regions $\{z:r_n\le|z|\le r_{n+1},~\theta_n<\arg z\le\hat\theta\}$, this showing that the sectors $\Sigma_\nu$ are the natural pole-free sectors in the sense that for every $\delta>0$, $S_\delta=\{z:|\arg z-\nu{\ts\frac\pi 2}|<{\ts\frac\pi 4}-\delta\}$ contains only finitely many poles.

For  $z\in S_\delta$ sufficiently large we set $\max\{|p(z)|,|q(z)|\}=\epsilon|z|$. Since
\be{p+q}|p(z)|+|q(z)|=o(|z|)\quad(z\to\infty {\rm~on~}S_{\frac\delta 2}),\ee
we may assume $\epsilon=\epsilon(z)<\frac12$. From (\ref{Hsystem}) and $|p'(z)|+|q'(z)|\to 0$
as $z\to\infty$ on $S_\delta$ (this following from (\ref{p+q}) and Cauchy's Theorem), hence $|p'(z)|+|q'(z)|<1$, say, we obtain
$$|zp(z)|< \epsilon^2|z|^2+|\alpha|+1\und |zq(z)|<\epsilon^2|z|^2+|\beta|+1.$$
This yields $(\epsilon-\epsilon^2)|z|^2< K$ and $\epsilon<2K|z|^{-2}$,
hence
$$|p(z)|+|q(z)|=O(|z|^{-1})\und |p'(z)|+|q'(z)|=O(|z|^{-2})\quad(z\to\infty {\rm ~on~}S_\delta),$$
again by Cauchy's Theorem. From (\ref{Hsystem}) it then follows that
$$zp(z)+\alpha=O(|z|^{-2})\und zq(z)+\beta=O(|z|^{-2})$$
as $z\to\infty$ on $S_\delta$. Now assume that
\be{ASYMP}\begin{array}{rcl}
p(z)&=&\ds-\alpha z^{-1}+\sum_{\nu=1}^na_\nu z^{-2\nu-1}+O(|z|^{-2n-3})=\phi_n(z)+O(|z|^{-2n-3})\cr
q(z)&=&\ds-\beta z^{-1}+\sum_{\nu=1}^nb_\nu z^{-2\nu-1}+O(|z|^{-2n-3})=\psi_n(z)+O(|z|^{-2n-3})
\end{array}\ee
has already been proved. From (\ref{Hsystem}) it then follows that
$$\begin{array}{rcl}
zp(z)+\alpha&=&\ds-q(z)^2-p'(z)=-\psi_n(z)^2-\phi_n'(z)+O(|z|^{-2n-4})\cr
zq(z)+\beta&=&\ds -p(z)^2+q'(z)=-\phi_n(z)^2+\psi_n'(z)+O(|z|^{-2n-4})
\end{array}$$
holds. Since $\phi_n'$, $\psi_n'$, $\phi_n^2$, and $\psi_n^2$ are even functions, (\ref{ASYMP}) holds with $n$ replaced by $n+1$.
A more detailed computation(\footnote{Many of the computations at various places were performed with the help of {\sf maple}.}) then gives~(\ref{ASYMP0}). \Ende

\subsection{Asymptotics on pole-free sectors with $\mathbf{\CL(p,q)=\{{\ts\frac13}\}}$} The other cases are yet easier to deal with, since the principal terms $-\tau z$ and $-\bar\tau z$ ($\tau^3=1$) are already known. Details are left to the reader. There is, however,
one remarkable difference: the natural pole-free sectors now are $\Sigma_\nu: \nu{\ts\frac\pi 2}<\arg z<(\nu+1){\ts\frac\pi 2}$, and these sectors are
bordered by sequences of poles $\lambda_n$ satisfying
(\ref{lambdarek}) with $\varpi=\ts\frac{2\pi}{\sqrt 3}$ and having counting function
$n(r,\{\lambda_n\})\sim \frac{\sqrt 3 r^2}{4\pi}.$
We also note that in the present case the algebraic curve (\ref{algcurve}) is reducible ($c={\ts\frac13})$.

\begin{theo}\label{ASY1}Suppose that $p(z)=-\tau z+o(|z|)$ and $q(z)=-\bar\tau z+o(|z|)$ $(\tau^3=1)$ hold as $z\to\infty$ on some ray
$\hat\sigma:\arg z=\hat\theta\not\equiv 0~\mod {\ts\frac\pi 2}$. Then $p$, $q$, and $H$ have asymptotic expansions
\be{ASYMP1}
\begin{array}{rcl}
p(z)&\sim&-\tau z+\frac{\alpha+2\bar\tau\beta+\tau}{3z}+\frac{3\alpha+2\bar\tau\beta+2\bar\tau\alpha^2+3\beta^2+4\tau\alpha\beta-2\tau}{9z^3}+\cdots\cr
q(z)&\sim&-\bar\tau z+\frac{\beta+2\tau\alpha-\bar\tau}{3z}-\frac{3\beta+2\tau\alpha-2\tau\beta^2-3\alpha^2-4\bar\tau\alpha\beta+2\bar\tau}{9z^3}+\cdots\cr
 H(z)&\sim&\frac{z^3}3-(\tau\alpha+\bar\tau\beta)z+\frac{\tau\alpha^2+\bar\tau\beta^2+\alpha\beta-1}{3z}+\cdots\end{array}\ee
as $z\to\infty$ on the sector $\Sigma_\nu: \nu{\ts\frac\pi 2}<\arg z<(\nu+1){\ts\frac\pi 2}$ that contains $\hat\sigma$.\end{theo}

\subsection{Asymptotics on adjacent pole-free sectors}Both theorems may be completed as follows:

\begin{theo}\label{ASY2}Suppose that the same asymptotics (\ref{ASYMP0}) and  (\ref{ASYMP1}) hold on adjacent sectors $\Sigma_\nu$ and $\Sigma_{\nu+1}$,
respectively. Then this is true on $\big(\overline{\Sigma_\nu\cup\Sigma_{\nu+1}}\big)^\circ.$\end{theo}

\proof In both cases, $g(z)=\exp \int H(z)\,dz$ is an entire function of finite order ($\le 4$) having simple zeros at the poles of $p$ and $q$. Taking into account the asymptotics (\ref{ASYMP0}) and (\ref{ASYMP1}) of $H$,  we consider
$$f(z)=z^{\alpha\beta}\,g(z)\und f(z)=g(z)\,\ts z^{-(\tau\alpha^2+\bar\tau\beta^2+\alpha\beta-1)/3}\,\exp(-\frac{z^4}{12}+\frac12(\tau\alpha+\bar\tau\beta)z^2)$$
in the respective cases. Then on $|\arg z-\bar\theta_\nu|<\frac\pi 8$, say ($\arg z=\bar\theta_\nu$ denotes the ray that separates $\Sigma_\nu$ and $\Sigma_{\nu+1}$), $f$ has order of growth at most four,
$$\limsup_{z\to\infty}\frac{\log\log|f(z)|}{\log|z|}\le 4\quad(|\arg z-\bar\theta_\nu|<\ts\frac\pi 8),$$
and satisfies $\lim\limits_{r\to\infty}f(re^{i(\bar\theta_\nu\pm\delta)})=C_\pm$
for $\delta=\frac\pi{10}$, say. The Phragm\'en-Lindel\"of Principle(\footnote{The version we will use can easily be derived from the following standard version:
\it Let $f$ be holomorphic on the domain $D:|z|>r_0,~\Im z>0$, with continuous boundary values on $(-\infty,-r_0)\cup(r_0,+\infty)$. If $f$ satisfies
$\ds\limsup\limits_{z\to\infty}\frac{\log|f(z)|}{|z|}=0$ on $D$, $\lim\limits_{x\to +\infty} f(x)=a$ and $\lim\limits_{x\to -\infty} f(x)=b$, then $a=b$ and $\lim\limits_{z\to \infty} f(z)=a$ holds on $D$.})
then shows that $C_+=C_-=C$, and $f(z)=C+o(1)$ holds as $z\to\infty$ on $|\arg z-\bar\theta_\nu|<\delta$. Hence $f$ and also $g$ has only finitely many zeros,
and $p$ has only finitely many poles on that sector, that is, $p$ and $q$ have an asymptotic expansion on $|\arg z-\hat\theta|<\delta$ which coincides with those on the sectors
$\Sigma_\nu$ and $\Sigma_{\nu+1}.$ \Ende

\subsection{B\"acklund transformations and asymptotics}Under the B\"acklund transformation ${\sf B}_\omega$ the asymptotics (\ref{ASYMP1}) changes as follows:
\be{ASYMPCHANGE}
p(z)\sim-\tau z {\rm ~and~}\left\{\begin{array}{c}\omega=\bar\tau\cr\omega\ne\bar\tau\end{array}\right. {\rm ~implies~}
{\sf B}_\omega p(z)\sim\left\{\begin{array}{c}-\tau z\cr-\bar\tau z.\end{array}\right.\ee

\subsection{Existence and uniqueness of asymptotic expansions}In Theorem~\ref{ASY0} and \ref{ASY1} it was shown that specific solutions have asymptotic expansions
on certain sectors of central angle ${\ts\frac\pi 2}$. This should not be confused with our next result about existence and uniqueness of solutions having asymptotic expansions.

\begin{theo}\label{EXUNIQUE} To any half-plane $|\arg z-(2\nu+1){\ts\frac\pi 4}|<{\ts\frac\pi 2}$ and $|\arg z-\nu{\ts\frac\pi 2}|<{\ts\frac\pi 2}$ there exists a unique solution to (\ref{Hsystem}) with prescribed asymptotics (\ref{ASYMP0}) and (\ref{ASYMP1}), respectively.\end{theo}

\proof To prove {\it existence} we  set, in the first case (\ref{ASYMP0}), $t=z^2$,
$p(z)=t^{-\frac12}u(t)^2$, and $q(z)=t^{-\frac12}v(t)^2$ to obtain
$$\dot u=-\frac{\alpha+u^2}{4u}+\frac{u^2-v^4}{4tu},~\dot v=\frac{\beta+v^2}{4v}+\frac{u^4+v^2}{4tv}.$$
Then given any half plane $H$, Theorem 14.1~in Wasow's monograph \cite{Wasow1} applies to the corresponding system for $x=u-\sqrt{-\alpha}$ and $y=v-\sqrt{-\beta}$ in $H$; it yields existence of some solution
having asymptotic expansion $u\sim\sqrt{-\alpha}+\sum\limits_{k=1}^\infty a_kt^{-k}$, $v\sim\sqrt{-\beta}+\sum\limits_{k=1}^\infty b_kt^{-k}$ as $t\to\infty$ on $H$.  Thus (\ref{Hsystem})
has some solution with asymptotic expansion (necessarily given by (\ref{ASYMP0})) on $S=\sqrt{H}$. If, however, $S$ contains some ray $\hat\sigma_\nu:\arg z=\nu{\ts\frac\pi 2}$, then by Theorem~\ref{ASY0}, the asymptotic expansion extends to the half-plane $|\arg z-\nu{\ts\frac\pi 2}|<{\ts\frac\pi 2}$.
We note, however, that the proof only works if $\alpha\beta\ne 0$. In the trivial case $\alpha=\beta=0$ we have $p\equiv q\equiv 0$.
In case of $\alpha=0$ and $\beta\ne 0$, say, we set $t=z^2$, $p(z)=t^{-3/2}u(t)^2$, and $q(z)=t^{-1/2}v(t)^2$ to obtain
the system
$$\dot u=-\frac{u^2+v^4}{4u}+\frac{3u}{4t},~\dot v=\frac{v^2+\beta}{4v}+\frac{v}{4t}+\frac{u}{4t^3v}.$$
Again, Theorem 14.1 in \cite{Wasow1} applies to the system for $x=u-i\beta$, $y=v-\sqrt{-\beta}$ ($u^2+v^4=v^2+\beta=0$ is solved by
$v=\sqrt{-\beta}$, $u=\pm i\beta$).

\smallskip In the same way one may deal with the second case (\ref{ASYMP1}): just set  $t=z^2$,
$p(z)=t^{\frac12}u(t)^2$, and $q(z)=t^{\frac12}v(t)^2$ to obtain
$$\dot u=-\frac{u^2+v^4}{4u}-\frac{\alpha+u^2}{4tu},~\dot v=\frac{u^4+v^2}{4v}+\frac{\beta-v^2}{4tv},$$
and apply  Theorem 14.1 in \cite{Wasow1} to the system for $x=u-\sqrt{-\tau}$, $y=v-\sqrt{-\bar\tau}$ (the non-trivial solutions to $\xi+\eta^2=\xi^2+\eta=0$ are $(\xi,\eta)=(-\tau,-\bar\tau)$ with $\tau^3=1$).

\medskip To prove {\it uniqueness} we consider two solutions $(p_1,q_1)$ and $(p_2,q_2)$ having the same asymptotic expansion on some sector $S$. Then
$u=p_1-p_2$ and $v=q_1-q_2$ tend to zero faster than every power $z^{-n}$, and $\u=\left(\!\!\!\begin{array}{c}u\cr v\end{array}\!\!\!\right)$ solves
$z^{-1}\u'=A(z)\u$ with
$A(\infty)=\left(\!\!\!\begin{array}{cc}-1\!&\!\!0\cr 0\!\!&\!\!1\end{array}\!\!\!\right)$ (eigenvalues $\pm 1$) and $A(\infty)=\left(\!\!\!\begin{array}{cc}-1\!\!&\!\!-2\bar\tau\cr -2\tau\!\!&\!\!1\end{array}\!\!\!\right)$ (eigenvalues $\pm \sqrt 5$),
respectively. In every sector $|\arg z-\nu{\ts\frac\pi 2}|<{\ts\frac\pi 4}$ and $|\arg z-(2\nu+1){\ts\frac\pi 4}|<{\ts\frac\pi 4}$, respectively, systems of this kind have some solution $\u_1$ that tends to infinity,
and also some non-trivial solution $\u_2$ (called sub-dominant) that tends to zero exponentially as $z\to\infty$.
If, however, $S$ is any sector with central angle greater than ${\ts\frac\pi 2}$, then {\it every} non-trivial solution to $z^{-1}\u'=A(z)\u$ tends to infinity as $z\to\infty$ on some sub-sector of $S$. This proves $(p_1,q_1)=(p_2,q_2)$ if $S$ is a half-plane. \Ende

\medskip{\sc Remark.} Theorem~14.1 in \cite{Wasow1} applies, in particular, to algebraic systems
$$t^{-n}\dot\y=A\y+\f(t,\y)$$
$(\det A\ne 0$, $\|\f(t,\y)\|=O(|t|^{-1})+O(\|\y\|^2)$ as $t\to +\infty,$ $\|\y\|\to 0)$
on arbitrary sectors with central angle $\frac\pi{n+1}$; the number $n+1\ge 1$ is called  {\it rank}. We note that it is also part of the hypothesis that the system has a {\it formal} solution $\sum\limits_{k=1}^\infty \c_kt^{-k}$, giving rise to the {\it asymptotic} solution $\y\sim\sum\limits_{k=1}^\infty \c_kt^{-k}$; existence of the formal solution, however, follows immediately from $\det A\ne 0$.

\subsection{Rational solutions}For rational solutions, the asymptotic series in (\ref{ASYMP0}) and (\ref{ASYMP1}) converge on $|z|>r_0$. Suppose that $p$ has $n^\varrho$ finite poles with residue $\varrho$. In the first case, integrating $p$, $q$, and $H$ over any large circle $|z|=r$ yields
$$\ts\sum\limits_{\varrho^3=1}\varrho n^\varrho=-\alpha,~\sum\limits_{\varrho^3=1}-\bar\varrho n^\varrho=-\beta\und\sum\limits_{\varrho^3=1}n^\varrho=-\alpha\beta,$$
hence $\beta=-\bar\alpha$, $\Re\alpha\in\frac12\Z$, $\Im\alpha\in \frac{\sqrt 3}2\Z$, and $-\alpha\beta=|\alpha|^2\in\N_0.$ In the second case (\ref{ASYMP1}) we obtain in the same way, restricting ourselves to the case $\tau=1$ in (\ref{ASYMP1}),
$$\ts\sum\limits_{\varrho^3=1}\varrho n^\varrho={\ts\frac13}(\alpha+2\beta+1),~\sum\limits_{\varrho^3=1}-\bar\varrho n^\varrho={\ts\frac13}(2\alpha+\beta-1),~\sum\limits_{\varrho^3=1}n^\varrho={\ts\frac13}(\alpha^2+\beta^2+\alpha\beta-1),$$
hence again $\beta=-\bar\alpha$, $\Re\alpha\in 1+\frac32\Z$, $\Im\alpha\in \frac{\sqrt 3}2\Z$, and $\alpha^2+\bar\alpha^2-|\alpha|^2-1\in 3\N_0$.

\section{Solutions of minimal growth}

\subsection{The class $\mathbf{\P_{\alpha,\beta}}$} Our focus is on transcendental solutions $(p,q)$ to (\ref{Hsystem}) satisfying $n(r,\Lambda)\asymp r^2$. These solutions form the class $\P_{\alpha,\beta}$.
If it is clear which parameters are in question or if the parameters are irrelevant, we shall write just $\P$.
By Theorem~\ref{EQUALN}, any solution that misses one kind of poles belongs to $\P$, while from Theorem~\ref{ORDERDICHO} it follows that either $\CL(p,q)=\{{\ts\frac13}\}$ or else $\CL(p,q)=\{0\}$ holds for $(p,q)\in\P$. We also note that the union of classes $\P_{\alpha,\beta}$ is invariant under B\"acklund transformations: $T(r,{\sf B}p)=O(T(r,p))=O(r^2)$.

\subsection{The case $\mathbf{\CL(p,q)=\{0\}}$}In the first case we have

\begin{theo}\label{ZWEI}The solutions $(p,q)$ to equation (\ref{Hsystem}) with cluster set $\CL(p,q)=\{0\}$
and counting function $n(r,\Lambda)=O(r^2)$ are rational functions, $g(z)=\exp\big(\int H(z)\,dz\big)$ is a polynomial and $p$ and $q$ are rational of degree $-\alpha\beta$.\end{theo}

\proof There exist only finitely many strings of poles, and $p$ and $q$ have asymptotic expansions given by (\ref{ASYMP0}), one and the same on each sector $|\arg z-\nu{\ts\frac\pi 2}|<{\ts\frac\pi 4}$.
From Theorem~\ref{ASY2} it then follows that these expansions hold throughout the whole plane, hence $p$ and $q$ are rational functions, while the entire function $g$ satisfies $g(z)\sim Cz^{-\alpha\beta}$ as $z\to\infty$, hence is a polynomial of degree $-\alpha\beta$. \Ende

\subsection{The case $\mathbf{\CL(p,q)=\{\frac13\}}$}For $(p,q)\in\P_{\alpha,\beta}$ with $\CL(p,q)=\{\frac13\}$, $p$ has $\n^\varrho(p)$ strings
of poles with residue $\varrho$, $\varrho^3=1$. It is quite natural to consider subclasses
of $\P_{\alpha,\beta}$ as follows: the pair $(p,q)$ belongs to class\medskip

\begin{tabular}{rl}
$\P_{\alpha,\beta}(\varrho_1)$,& if $\n^{\varrho_1}(p)>0$, but $\n^{\varrho_2}(p)=\n^{\varrho_3}(p)=0$ ($\varrho_\nu^3=1$);\cr
$\P_{\alpha,\beta}(\varrho_1,\varrho_2)$,& if $\n^{\varrho_1}(p)>0$ and $\n^{\varrho_2}(p)>0$,
but $\n^{\varrho_3}(p)=0$.
\end{tabular}

\medskip It is obvious that the six classes $\P_{\alpha,\beta}(\varrho)$ and $\P_{\alpha,\beta}(\varrho_1,\varrho_2)$ are mutually disjoint, and are contained in $\P_{\alpha,\beta}$ by Theorem~\ref{EQUALN}. We also note that $\n^\varrho(p)=0$
implies that $p$ has at most finitely many poles with residue $\varrho$.

\subsection{Asymptotics and the distribution of residues}Let $(p,q)$ be any pair in $\P$ with cluster set $\CL(p,q)=\{\frac13\}$. With each string of poles $s=(\lambda_k)$ we associate a polygon $\pi(s)$ with vertices $\lambda_k$. We assume that the polygons $\pi_\nu=\pi(s_\nu)$ are in cyclic order (counter-clockwise).
They divide $|z|>r_0$ into finitely many domains $D_1, D_2, \ldots D_n$ ($D_{n+1}=D_1$); $D_\nu$ is bounded by $\pi_\nu$, $\pi_{\nu+1}$ and some sub-arc of $|z|=r_0$.
Re-scaling along any sequence $(\kappa_m)$ in $D_\nu$ with $|\kappa_m|\dist(\kappa_m,\Lambda)\to\infty$
yields constant limit functions $\p=-\omega_\nu$ and $\q=-\bar\omega_\nu$, hence
\be{asy1}p(z)=-\omega_\nu z+o(|z|),~ q(z)=-\bar\omega_\nu z+o(|z|)\quad(z\in D_\nu,~|z|\dist(z,\partial D_\nu)\to\infty);\ee
this does, of course, not mean that $p$ has an asymptotic expansion on $D_\nu$, except when $D_\nu$ contains some sector. To determine $\varrho_\nu=\Res_{s_\nu} p$ we assume $\arg \lambda_k\sim 0$ for the sake of simplicity, and compute
\be{RESINT}\frac1{2\pi i}\int_{\gamma_\nu} p(z)\,dz=\varrho_\nu\frac{\sqrt 3}{4\pi}r^2+o(r^2)\ee
along a closed curve $\gamma_\nu$ surrounding $s_\nu\cap\{z:r_0<|z|<r\}$. Given $\delta>0$ sufficiently small, $\gamma_\nu$ consists of sub-arcs $\sigma_{r_0}$ and $\sigma_r$ of $|z|=r_0$ and $|z|=r$, of length at most $\delta r_0$ and $\delta r$, respectively, and arcs $\gamma'_{\nu-1}$ and $\gamma'_\nu$ in $D_{\nu-1}$ and $D_\nu$, such that $p(z)=\omega_{\nu-1}z+o(|z|)$ and $p(z)=\omega_{\nu}z+o(|z|)$ hold as $z\to\infty$ on $\gamma'_{\nu-1}$  and $\gamma'_\nu$, respectively:

{\small $$\begin{array}{c}
p(z)=-\omega_{\nu}z+o(|z|)\cr
\leftarrow\gamma'_{\nu}\cr
\sigma_{r_0}\downarrow\begin{array}{|cccccccccc|}\hline
&&&&&&&&&\cr
\bullet&\bullet&\bullet&\bullet&\bullet&\bullet&\bullet&\bullet&\bullet&\bullet\cr
&&&&&&&&&\cr\hline\end{array}\uparrow\sigma_{r}\cr
\gamma'_{\nu-1}\rightarrow \cr
p(z)=-\omega_{\nu-1}z+o(|z|)\end{array}$$
\begin{center}{\sc Figure 1.} The poles $\bullet$ have residue $\varrho_\nu$.\end{center}}

These arcs contribute to the integral (\ref{RESINT}) as follows:
$$\frac1{2\pi i}\int_{\sigma_{r_0}}p(z)\,dz+\frac1{2\pi i}\int_{\sigma_{r}}p(z)\,dz =O(\delta r^2),$$
$$\frac1{2\pi i}\int_{\gamma'_{\nu-1}}p(z)\,dz+\frac1{2\pi i}\int_{\gamma'_{\nu}}p(z)\,dz=\frac{-\omega_{\nu-1}+\omega_\nu}{i\sqrt 3}\,\frac{\sqrt 3}{4\pi}r^2+O(\delta r^2),$$
which implies $\ds\varrho_\nu=i\frac{\omega_{\nu-1}-\omega_\nu}{\sqrt 3}$. Table~1 displays the various possibilities (yes/no).
$$\begin{array}{|c|c|c|c|}\hline
\omega_{\nu-1}&\omega_\nu&\varrho_\nu&?\cr\hline
1&\omega&-\bar\omega&{\rm no}\cr\hline
1&\bar\omega&\omega&{\rm yes}\cr\hline
\omega&\bar\omega&-1&{\rm no}\cr\hline
\omega&1&\bar\omega&{\rm yes}\cr\hline
\bar\omega&\omega&1&{\rm yes}\cr\hline
\bar\omega&1&-\omega&\rm no\cr\hline\end{array}\qquad\qquad
\begin{array}{|c|c|c|c|}\hline
\omega_{\nu-1}&\omega_\nu&\varrho_\nu&?\cr\hline
1&\omega&\bar\omega&{\rm yes}\cr\hline
1&\bar\omega&-\omega&{\rm no}\cr\hline
\omega&\bar\omega&1&{\rm yes}\cr\hline
\omega&1&-\bar\omega&{\rm no}\cr\hline
\bar\omega&\omega&-1&{\rm no}\cr\hline
\bar\omega&1&\omega&\rm yes\cr\hline\end{array}$$

\begin{center}\quad\small $\arg p_k\sim 0$ mod $\pi$\qquad\qquad\qquad $\arg p_k\sim {\ts\frac\pi 2}$ mod $\pi$\\
{\sc Table~1.} The asymptotics determines the residues, and vice versa ($\omega={\ts\frac12(-1+i\sqrt 3)}$).
\end{center}

\subsection{Solutions of the first kind}Suppose that $(p,q)\in\P_{\alpha,\beta}(1)$, say.
Then $w=p+q-z$ has only finitely many poles and vanishes at every pole with residue $1$.
On the other hand it follows from $m(r,w)=O(\log r)$ that $w$ is a rational function, hence vanishes identically. From equation (\ref{PIVmod}), (\ref{Hsystem}), and also from the development of $p+q-z$ in (\ref{POLESDiv}) it then easily follows that
$\alpha-\beta+1=0$ and $\h(\lambda)=0$. Conversely, let $(p,q)$ be any solution to (\ref{Hsystem}) such that $p(z_0)+q(z_0)-z_0=0$ at some regular point. Then $w=p+q-z$ satisfies
$$w'=-q^2+p^2-z(p-q)-\alpha+\beta-1=(p(z)-q(z))w\und w(z_0)=0,$$
hence vanishes identically ($p(z)-q(z)$ may be viewed as a coefficient).

\begin{theo}\label{FK}The class $\P_{\alpha,\beta}(\varrho)$ is non-empty if and only if $\bar\varrho\alpha-\varrho\beta+1=0$.
In that case, $(p,q)\in\P_{\alpha,\beta}(\varrho)$ satisfies
$$\bar\varrho p+\varrho q-z=0,\quad\h(\lambda)=0,\und\ts H=\frac13 z^3+\bar\varrho z+\bar\varrho p,$$
and $p$ solves the Riccati equation
\be{RICCp}p'=-\alpha-\bar\varrho z^2+z p-\bar\varrho p^2.\ee\end{theo}

\proof We restrict ourselves to the case $\varrho=1$. The differential equation (\ref{RICCp}) and the form of $H$ is obtained from (\ref{Hsystem}) and (\ref{hamiltonian}) by replacing $q$ by $z-p$.
To prove that $\P_{\alpha,\alpha+1}(1)$ is non-empty, we start with any transcendental solution to equation (\ref{RICCp}) -- not {\it every} solution can be rational --, and set $q=z-p$.
Then $p'=-q^2-zp-\alpha$ and $q'=p^2+zq+\alpha+1$ immediately follow. \Ende

\medskip{\sc Remark.} For $(p,q)\in\P_{\alpha,\varrho\alpha+\bar\varrho}(\varrho)$ ($\varrho=\frac12(-1+i\sqrt 3)$),
$w=p+q-z$ has simple poles with residue $\Res_\lambda w=i\sqrt 3$, and  $w'+\frac{w^2}{i\sqrt 3}-i\sqrt 3 w$ is regular at $\lambda$ and assumes the value
$\frac12((3-i\sqrt 3)\alpha-(3-i\sqrt 3))$ at $z=\lambda$. This implies
$$w'={\ts\frac12}((3-i\sqrt 3)\alpha-(3-i\sqrt 3))+i\sqrt 3 w-\frac{w^2}{i\sqrt 3},$$
and $w'=\frac{i}{\sqrt 3}w(w+3)$ if $\alpha=1$. In the latter case, $w$ has no zeros.

For $\varrho=1$, say, equation (\ref{RICCp}) may be transformed into
$$u'=\ts\pm i\frac{2\alpha+1}{\sqrt 3}+\zeta^2-u^2$$
(set $\zeta=az$ and $p(z)-\frac z2=au(\zeta)$ with $a^4=-\frac34$).
The transcendental solutions to equations of this type were analysed in \cite{NStRiccati}, two major results may be described as follows:
\begin{itemize}\item The set of poles of {\it generic} solutions $u$ consist of {\it four} sequences $(\zeta^{(\nu)}_n)$ such that $\arg\zeta^{(\nu)}_n\sim (2\nu-1){\ts\frac\pi 4}$
and $|\zeta^{(\nu)}_n|\sim \sqrt{2\pi n}$ as $n\to\infty$. Moreover, $u$ has an asymptotic expansion
$u(\zeta)=(-1)^{\nu}\zeta+O(|\zeta|^{-1})$ as $\zeta\to\infty$ on each open sector $|\arg \zeta-\nu{\ts\frac\pi 2}|<{\ts\frac\pi 4}$.
\item There are four {\it degenerate} solutions $u_\mu$ having only sequences of poles for $\nu=\mu$ and $\nu=\mu+1$, such that $u_\mu(\zeta)\sim (-1)^\mu\zeta$ holds on
$|\arg \zeta-\mu{\ts\frac\pi 2}|<{\ts\frac\pi 4}$, while $u_\mu(\zeta)\sim -(-1)^\mu\zeta$ holds on $2(\mu+1){\ts\frac\pi 4}<\arg\zeta<(2\mu+7){\ts\frac\pi 4}$
(for some parameters $\alpha$, rational  instead of degenerate solutions might occur).\end{itemize}

In our case (\ref{RICCp}) it follows that $\Lambda$ consist of two or four
strings $(\lambda_k)$, each with counting function $n(r,\{p_k\})\sim\frac{\sqrt 3}{4\pi} r^2$ and asymptotic to the rays $\arg z=\nu{\ts\frac\pi 2}$ ($\nu=\mu,\mu+1$ or $\nu=0,1,2,3$).
From Table~1 it follows that the residues determine the asymptotics.

{\small$$\begin{array}{ccc}
p\sim-\bar\omega z &\vertikal{\symbol}&p\sim-\omega z\cr
\horizontal{\symbol}\!\!\!\!\!\!\!&&\!\!\!\!\!\!\!\horizontal{\symbol}\cr
p\sim-\omega z   &\vertikal{\symbol}&p\sim-\bar\omega z \cr
\end{array}\qquad\qquad\begin{array}{ccc}
p\sim-\bar\omega z &\vertikal{\symbol}&p\sim-\omega z\cr
&&\!\!\!\!\!\!\!\horizontal{\symbol}\cr
p\sim-\bar\omega z&\phantom{\vertikal{\symbol}}&p\sim-\bar\omega z \cr
\end{array}$$

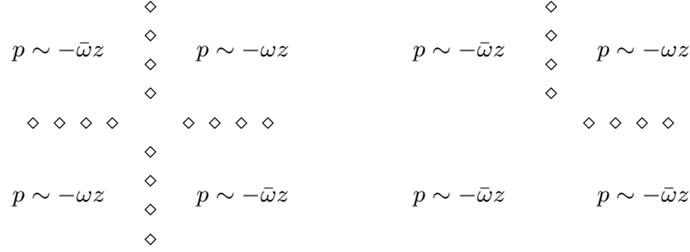
\begin{center}{\sc Figure~2.} Generic and degenerate distribution of poles and asymptotics in  $\P_{\alpha,\alpha+1}(1)$;\\
each symbol $\symbol$ represents a pole with residue $1$; $\omega={\ts\frac12(-1+i\sqrt 3)}$.\end{center}}

\subsection{Solutions of the second kind}We will confine ourselves to the symmetric class $\P_{\alpha,\beta}(\varrho,\bar\varrho)$, noting that ${\sf M}_{\bar\varrho}\P_{\alpha,\beta}(1,\varrho)=\P_{\varrho\alpha,\bar\varrho\beta}(\varrho,\bar\varrho)$. For $\omega\in\{\varrho,\bar\varrho\}$ we consider the linear combination $u_\omega=\omega p+\bar\omega q+\frac z2$; $u_\omega$ is regular at poles of $p$ with residue $\bar\omega$ and has simple poles at poles of $p$ with residue $\omega$:
\be{LKu}\begin{array}{rcl}
u_\omega(z)&=&\ds\frac{-2i\Im\omega}{z-\lambda}+a(\lambda;\omega)(z-\lambda)+\cdots\cr
a(\lambda;\omega)&=&-\ts\frac 12+\frac13[(\omega+2)\alpha-(\bar\omega+2)\beta]+\frac i2(\Im\omega)\lambda^2.
\end{array}\ee
Like in the first case we may conclude that $u_\omega$ satisfies the Riccati equation
\be{RICCu}u_\omega'=3a(z;\omega)+\frac{u_\omega^2}{2i\Im\omega}.\ee
This means, in particular, that $p$ has four [two] strings of poles with residues $\varrho$ and $\bar\varrho$, asymptotic to all [two adjacent] directions $\arg z=\nu\frac\pi2$.
Combining (\ref{Hsystem}) and (\ref{RICCu}) yields an algebraic equation $K_\omega(z,p,q)=0$, and the system $(K_\varrho=0,K_{\bar\varrho}=0)$ has resultant
$-\frac{16}3(\alpha-\beta-2)^2$ with respect to each of $p$ and $q$, thus
\be{BED2}-\alpha+\beta+2=0\ee
and also
\be{BEDK}K(z,p,q)=p^2+q^2+z^2-pq+z(p+q)+3\beta+3\equiv 0\ee
is necessary for $\P_{\alpha,\beta}(\varrho,\bar\varrho)\ne\emptyset$.
We also note that $(\alpha,\beta)=(1,-1)$ is the unique solution to the linear system
$-\alpha+\beta+2=0,$ $\bar\omega\alpha-\omega\beta+1=0.$
By Proposition~\ref{PROP1} below we obtain $(\alpha,\beta)\ne(1,-1)$ as further necessary condition for $\P_{\alpha,\beta}(\varrho,\bar\varrho)\ne\emptyset.$
Based on these considerations we will prove:

\begin{theo}The class $\P_{\alpha,\beta}(\varrho,\bar\varrho)$ with $(\alpha,\beta)\ne (1,-1)$ is non-empty if and only if (\ref{BED2}) holds.
The B\"acklund transformation ${\sf B}_\varrho$ maps $\P_{\beta+2,\beta}(\varrho,\bar\varrho)$ into
$\P_{\varrho\beta-\bar\varrho,\bar\varrho\beta+\bar\varrho-1}(\varrho)$. Moreover,
\be{BED2a}H(z,p,q)=q-p+{\ts\frac 13}z^3+(\beta+1)z\ee
and  $\h(\lambda)+i\Im\omega=(\beta+1)\lambda$ $(\Res_\lambda p=\omega\in\{\varrho,\bar\varrho\})$ hold.\end{theo}

\proof Replacing $pq$, $q^2+zp$, and $p^2+zq$ in (\ref{BEDK}) by $H'$, $-p'-2-\beta$, and $q'-\beta$, respectively, we obtain
$(q-p- H)'+z^2+\beta+1=0$, hence also
$$H(z,p,q)=q-p+{\ts\frac 13}z^3+(\beta+1)z+c$$
holds. Since $p$ and $q$ are transcendental, this is compatible with (\ref{BEDK}) if and only if $c=0$ (we note that
$\frac13(p^3+q^3-z^3)+zpq=\frac13(p+q-z)(p^2+q^2+z^2-pq+z(p+q)+2\beta+2)$).
Thus also (\ref{BED2a}) holds, while $\h(\lambda)+i\Im\omega=(\beta+1)\lambda$ at poles of $p$ with residue $\omega\in\{\varrho,\bar\varrho\}$ follows from (\ref{RICCu}) and (\ref{POLESH}).

To prove the second assertion, we note that $(\alpha,\beta)\ne(1,-1)$ and (\ref{BED2}) imply that $\varrho\alpha-\bar\varrho\beta+1\ne 0$, and also that the function $\varrho p+\bar\varrho q-z$ does not vanish identically since it has poles at poles of $p$ with residue $\varrho$.
According to (\ref{CHANGEofRESIDUES}), the poles of $p$ with residue $\varrho$ survive as poles of ${\sf B}_\varrho p$ with the same residue,
while the poles of $p$ with residue $\bar\varrho$ disappear (we note that $\n^\varrho(p)=\n^{\bar\varrho}(p)$). The zeros of $f=\varrho p+\bar\varrho q-z$ which are not poles of $p$  with residue $\bar\varrho$ may create poles with residue $\bar\varrho$.
From $T(r,f)=N(r,f)+O(\log r)=N(r,\Lambda^\varrho)+O(\log r)=N(r,\Lambda^{\bar\varrho})+o(r^2)$ and Nevanlinna's First Main Theorem, however, it follows that
$$N(r,\Lambda^{\bar\varrho})\le N(r,1/f)\le T(r,f)+O(1)=N(r,\Lambda^\varrho)+O(\log r)=N(r,\Lambda^{\bar\varrho})+o(r^2),$$
hence $f$ has at most $o(r^2)$ zeros on $|z|<r$ except those arising from poles of $p$ with residue $\bar\varrho$. This proves $\n^{\bar\varrho}({\sf B}_{\varrho}p)=\n^{1}({\sf B}_{\varrho}p)=0$ and ${\sf B}_{\varrho}(p,q)\in\P_{\varrho\beta-\bar\varrho,\bar\varrho\beta+\bar\varrho-1}(\varrho)$.

To prove that (\ref{BED2}) together with $(\alpha,\beta)\ne(1,-1)$ is sufficient for $\P_{\alpha,\beta}(\varrho,\bar\varrho)\ne\emptyset$, we start with any transcendental pair $(\tilde p,\tilde q)\in\P_{\tilde\alpha,\tilde\beta}(\varrho)$ with $(\tilde\alpha,\tilde\beta)=(\varrho\beta-\bar\varrho,\bar\varrho\beta+\bar\varrho-1)$, and note that $\bar\varrho\tilde\alpha-\varrho\tilde\beta+1=\beta-\varrho-(\beta+1-\varrho)+1=0$,
but $\varrho\tilde\alpha-\bar\varrho\tilde\beta+1=(\bar\varrho-\varrho)\beta-1-(\varrho-\bar\varrho)+1\ne 0$ for $\beta\ne -1$.
Thus $\bar\varrho \tilde p+\varrho p-z$ vanishes identically, but  $\varrho \tilde p+\bar\varrho p-z$ does not, so that $(p,q)={\sf B}_\varrho(\tilde p,\tilde q)$ is well defined.
From (\ref{CHANGEofRESIDUES}) it then follows that the poles of $\tilde p$ remain poles of $p$
with residue $\varrho$, while new poles with residue $\bar\varrho$, but none with residue $1$, might be created by the zeros of $\tilde f=\varrho \tilde p+\bar\varrho p-z$, if any. Since $(p,q)$ solves (\ref{Hsystem}) with parameter $(\beta+2,\beta)$, hence $(\beta+2)-\beta+1\ne 0$, the poles with residue $\bar\varrho$ actually exist, and $(p,q)\in\P_{\beta+2,\beta}(\varrho,\bar\varrho)$ follows. \Ende

{\small$$\begin{array}{ccc}
p\sim-\varrho z &\vertikal{\circ}\!\!\!\vertikal{\bullet}&p\sim-\bar\varrho z\cr
 \horizontal{\circ}\!\!\!\!\!\!&&\!\!\!\!\!\!\horizontal{\bullet}\cr
 \horizontal{\bullet}\!\!\!\!\!\!&&\!\!\!\!\!\!\horizontal{\circ}\cr
p\sim-\bar\varrho z  &\vertikal{\bullet}\!\!\!\vertikal{\circ}&p\sim-\varrho z \cr
 \end{array}$$
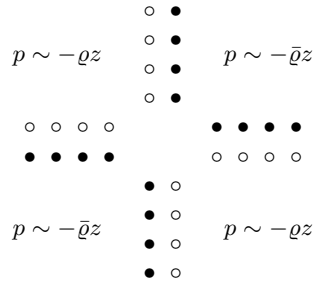
\begin{center}{\sc Figure~3.} Distribution of poles and asymptotics in  $\P(\varrho,\bar\varrho)$,
$\Res_\bullet p=\varrho=\frac12(-1+i\sqrt 3)$, $\Res_\circ p=\bar\varrho$; in the narrow domains between two strings, $p(z)=-z+o(|z|)$ holds.\end{center}}

\medskip{\sc Remark.} The case $(\alpha,\beta)=(1,-1)$ is indeed exceptional: $\P_{1,-1}(\varrho,\bar\varrho)$ is empty as follows from Proposition~\ref{PROP1} below with $\varrho_1=\varrho$ and $\varrho_1=\bar\varrho$, respectively. We also note that in this case, (\ref{Hsystem}) admits polynomial solutions $p=q=-z$.

\section{Solutions of the third kind}
\subsection{Preliminary remarks}Functions $(p,q)\in\P_{\alpha,\beta}$ have cluster set $\CL(p,q)=\{{\ts\frac13}\}$ and finitely many strings of poles, each  asymptotic to some ray $\arg z=\nu{\ts\frac\pi 2}$ and with counting function $\sim\frac{\sqrt 3}{4\pi}r^2$, so that
$$n(r,\Lambda)\sim \n(p)\frac{\sqrt 3}{4\pi}r^2,\quad n(r,\Lambda^\varrho)\sim \n^\varrho(p)\frac{\sqrt 3}{4\pi}r^2\und\sum_{\varrho^3=1}\n^\varrho(p)=\n(p).$$
Moreover, $p$ and $q$ have asymptotic expansions $p(z)\sim -\tau_\nu z$ and $q(z)\sim -\bar\tau_\nu z$ on the quarter-planes $(\nu-1){\ts\frac\pi 2}<\arg z<\nu{\ts\frac\pi 2}$ ($1\le\nu\le 4)$.
This will be abbreviated by writing
\be{SCHEMA}p:~\schema{\tau_1}{\tau_2}{\tau_3}{\tau_4}.\ee

\subsection{Special cases}Before proceeding further we will examine two special cases:

\begin{prop}\label{PROP1}Suppose that $(p,q)\in\P_{\alpha,\beta}$ satisfies
$$\bar\varrho_1\alpha-\varrho_1\beta+1=0\und\n^{\varrho_1}(p)>\ts\frac12(\n^{\varrho_2}(p)+\n^{\varrho_3}(p)).$$
Then $(p,q)\in\P_{\alpha,\beta}(\varrho_1)$, that is,  $\bar\varrho_1 p+\varrho_1 q-z\equiv 0$.\end{prop}

\proof We may assume that $\varrho_1=1$, hence $\alpha-\beta+1=0$ and $\n^1(p)>\frac12(\n^\varrho(p)+\n^{\bar\varrho}(p))$ with $\varrho=\frac12(-1+\sqrt 3)$
hold. Then $w=p+q-z$ vanishes at least twice at poles $\lambda\in\Lambda^1$, and
if we assume that $w$ does not vanish identically, it follows from $m(r,w)=O(\log r)$ and Nevanlinna's First Main Theorem that
$$2N(r,\Lambda^1)\le N(r,1/w)\le T(r,w)+O(1)=N(r,\Lambda^\varrho)+N(r,\Lambda^{\bar\varrho})+O(\log r),$$
hence also $2\n^1(p)\le \n^\varrho(p)+\n^{\bar\varrho}(p)$ against our hypothesis, and Proposition \ref{PROP1} is proved completely. \Ende

\begin{prop}\label{PROP2}For $(p,q)\in\P_{\alpha,\beta}$ with signature (\ref{SCHEMA}), the sum $\sum\limits_{\nu=1}^4(-1)^\nu\tau_\nu$ is non-zero.\end{prop}

\proof Obviously, $\sum\limits_{\nu=1}^4(-1)^\nu\tau_\nu=0$ is only possible if, up to some rotation of the plane and up to replacing $(p,q)$ by ${\sf M}_\varrho(p,q)$
with $\varrho$ suitably chosen, $(p,q)$ has the signature $\schema\tau\tau{\bar\tau}{\bar\tau}$ with $\tau={\ts\frac12(-1+i\sqrt 3)}$. From Theorem~\ref{ASY0} it then follows that
$$\ds H(z)-\frac{z^3}3\sim\left\{\begin{array}{ll}
\ds-(\tau\alpha+\bar\tau\beta)z+\frac{\tau\alpha^2+\bar\tau\beta^2+\alpha\beta-1}{3z}+\cdots&(\Im z>0)\cr
\ds-(\bar\tau\alpha+\tau\beta)z+\frac{\bar\tau\alpha^2+\tau\beta^2+\alpha\beta-1}{3z}+\cdots&(\Im z<0)\end{array}\right.$$
holds. The entire function $g(z)=\ts\exp\big(\int  H(z)\,dz-\frac{z^4}{12}\big)$
has order of growth at most four and satisfies
$$g(z)=\left\{\begin{array}{ll}Ae^{az^2}z^\mu(1+o(1))&(z=re^{i\delta}\to\infty)\cr
Be^{bz^2}z^\nu(1+o(1))&(z=re^{-i\delta}\to\infty)\end{array}\right.~(0<\delta<\pi),$$
with $AB\ne 0$, $a=-\frac12(\tau\alpha+\bar\tau\beta)$, $b=-\frac12(\bar\tau\alpha+\tau\beta)$, $\mu=\frac13(\tau\alpha^2+\bar\tau\beta^2+\alpha\beta-1)$ and $\nu=\frac13(\bar\tau\alpha^2+\tau\beta^2+\alpha\beta-1)$.

\medskip In the first step we will show that $\Re a=\Re b$, and assume to the contrary that $\Re a>\Re b$. Then
$f(z)=e^{-az^2}z^{-\mu}g(z)$ satisfies
$$f(re^{i\delta})\to A\ne 0\und |f(re^{-i\delta})|=O(e^{(\Re b-\Re a)r^2\cos 2\delta }r^M)\to 0$$
as $r\to\infty$, which for $0<\delta<\delta_0$ sufficiently small contradicts the
Phragm\'en-Lindel\"of Principle.
Thus $\Re a\le \Re b$ holds, and in the same manner (consider $\overline{f(\bar z)}$ instead of $f$ with $(a,b)$ replaced by $(\bar b,\bar a$)
we obtain $\Re b\le\Re a$. The same argument will now show that $\Im a\le\Im b$: assuming to the contrary that ($\Re a=\Re b$) and $\Im a>\Im b$ holds, it follows that
$$f(re^{i\delta})\to A\ne 0\und |f(re^{-i\delta})|=O(e^{-(\Im a-\Im b)r^2\sin 2\delta }r^M)\to 0\quad(r\to\infty),$$
which for $\delta>0$ sufficiently small again contradicts the Phragm\'en-Lindel\"of Principle.
 Moreover, $\Im a=\Im b$ implies
that $f$ has only finitely many zeros, and $p$ has only finitely many poles asymptotic to the positive real axis.
The same argument, however, works for $g(-z)$ on $|\arg z|<\delta_0$, with $(a,b)$ replaced by $(b,a)$, this showing that also $\Im b\le\Im a$ holds. In particular, we have
$a=b$, this implying that $g$ has only finitely many zeros at all, and $p$ has only finitely many poles in contrast to our assumption that $p$ is transcendental.
Thus Proposition \ref{PROP2} is completely proved. \Ende

\subsection{The class $\mathbf{\P_{\alpha,\beta}}$.} We will now prove the main result of this section:
\begin{theo}\label{TK}To any $(p,q)\in\P_{\alpha,\beta}$ there exists some B\"acklund transformation {\sf B}, such that
$${\sf B}p+{\sf B}q-z=0\und {\sf B}\alpha-{\sf B}\beta+1=0$$
holds, that is, ${\sf B}(p,q)\in\P_{{\sf B}\alpha,{\sf B}\alpha+1}(1)$.\end{theo}

\proof Suppose $(p,q)\in\P_{\alpha,\beta}$ has signature (\ref{SCHEMA}). On combination with the asymptotics of $p$, Cauchy's Theorem gives
$$\sum\limits_{\varrho^3=1}\varrho\n^\varrho(p)=2i\sum\limits_{\nu=1}^4(-1)^\nu\frac{\tau_\nu}{\sqrt 3}.$$
Taking real and imaginary parts, we obtain
\be{REALIMAG}\begin{array}{rcl}
\n^1(p)-\frac12(\n^\varrho(p)+\n^{\bar\varrho}(p))&=&-2\ds\sum\limits_{\nu=1}^4(-1)^\nu\frac{\Im\tau_\nu}{\sqrt 3}\cr
\frac{1}2(\n^\varrho(p)-\n^{\bar\varrho}(p))&=&\phantom{-}2\ds\sum\limits_{\nu=1}^4(-1)^\nu\frac{\Re\tau_\nu}3\end{array}\quad(\varrho={\ts\frac12(-1+i\sqrt 3)}).\ee
By Proposition~\ref{PROP2} we have
$\n^{\varrho_1}(p)>\frac 12(\n^{\varrho_2}(p)+\n^{\varrho_3}(p))$ for some $\varrho_1$. We may assume $\varrho_1=1$, hence
\be{N1GROSS}\n^1(p)>\ts\frac12(\n^\varrho(p)+\n^{\bar\varrho}(p))>0;\ee
otherwise we would replace $(p,q)$ by ${\sf M}_{\varrho_1}(p,q)=(\bar\varrho_1 p,\varrho_1 q)$,
but retain the notation $(p,q)$ and $(\alpha,\beta)$. By Proposition~\ref{PROP1} we are done if $\alpha-\beta+1=0$.
Otherwise we apply ${\sf B}_1$ to obtain
$$\tilde p(z)=\ds p(z)-\frac{\alpha-\beta+1}{p(z)+q(z)-z},~\tilde q(z)=q(z)+\frac{\alpha-\beta+1}{p(z)+ q(z)-z}.$$
Under ${\sf B}_1$, the Hamiltonian and the signature change as follows:
$$\tilde H(z)- H(z)=-\frac{\alpha-\beta+1}{ p(z)+ q(z)-z}=\tilde p(z)-p(z)=-\tilde q(z)+q(z)$$
and
$$p:~ \schema{\tau_1}{\tau_2}{\tau_3}{\tau_4}~ {\buildrel{\sf B}_1\over\longrightarrow}~  \tilde p:~ \schema{\bar\tau_1}{\bar\tau_2}{\bar\tau_3}{\bar\tau_4}$$
Although $H$ and $\tilde H$ cannot be controlled sufficiently well close to the rays $\arg z=\nu{\ts\frac\pi 2}$, this can be done for $\tilde H(z)- H(z)$ to compute
$$n(r,\tilde p)- n(r,p)=\frac 1{2\pi i}\int_{\Gamma_r}( \tilde H(z)- H(z))\,dz=\frac 1{2\pi i}\int_{\Gamma_r}(\tilde p(z)-p(z))\,dz$$
asymptotically (it is obvious that the curve $\Gamma_r$, constructed in section \ref{VD}, may be modified in such a way that  $|p(z)|+|\tilde p(z)|=O(|z|)$ holds on $\Gamma_r$).
From
$$\tilde p(z)-p(z)=2i(\Im\tau_\nu)z+O(|z|^{-1}),$$
uniformly on $|\arg z-(\nu-\frac 12)\pi|<{\ts\frac\pi 4}-\delta,$ and (\ref{N1GROSS}) it follows that
\be{CHANGE}\n(\tilde p)-\n(p)=4\sum\limits_{\nu=1}^4(-1)^{\nu}\frac{\Im\tau_\nu}{\sqrt 3}=-2\n^1(p)+\n^\varrho(p)+\n^{\bar\varrho}(p)\le -1.\ee
Repeating this argument yields a sequence $({\sf B}^{[k]})$ of B\"acklund transformations ${\sf B}^{[k]}={\sf B}_1{\sf M}_{\varrho^{[k-1]}}{\sf B}^{[k-1]}$,
such that, for some $n$,
$${\sf B}^{[n]}p+{\sf B}^{[n]}q-z\equiv 0\und{\sf B}^{[n]}\alpha-{\sf B}^{[n]}\beta+1=0.~\ende$$

\subsection{Concluding remark}
Since the B\"acklund transformations do not change the number of different $\tau_\nu$'s, and the functions in $\P_{\alpha,\beta}(\omega)$ have only signatures with two
different $\tau_\nu$'s, this is also true for any $(p,q)\in\P$, and will simplify matter considerably. Up to a rotation of the plane, there are two possibilities: $\schema{\tau_1}{\tau_2}{\tau_1}{\tau_2}$ (generic) and $\schema{\tau_1}{\tau_2}{\tau_2}{\tau_2}$ (degenerate).
We will discuss in detail the case
$\n^1\le\n^{\bar\varrho}\le\n^\varrho$ $(\ts\varrho=\frac12(-1+i\sqrt 3)),$
which by (\ref{REALIMAG}) is equivalent to $\Im\tau_1\le \Im\tau_2$ and $\Re\tau_1\le\Re\tau_2$, but $\tau_1\ne\tau_2$.
The case $\n^1=\n^\varrho=\n^{\bar\varrho}$ is excluded by Proposition~\ref{PROP2}. There are two (non-degenerate) kinds of solutions:
\be{n1klein}{\rm (A)}~\schema{\bar\varrho}1{\bar\varrho}1 {\rm ~with~}\n^1=\n^{\bar\varrho}=\n^\varrho-4\und
{\rm (B)}~\schema{\bar\varrho}\varrho{\bar\varrho}\varrho {\rm ~with~}\n^\varrho=\n^{\bar\varrho}=\n^1+4.\ee
Based on Table 1, it is easy to show that, in first non-trivial cases, the residues are distributed as follows:
{\small$$
\begin{array}{ccc}
p\sim -z &\vertikal{\bullet}\!\!\!\vertikal{\symbol}\!\!\!\vertikal{\circ}\!\!\!\vertikal{\bullet}&p\sim -\bar\varrho z\cr
 \horizontal{\bullet}\!\!\!\!\!\!& &\!\!\!\!\!\!\horizontal{\bullet}\cr
    \horizontal{\symbol}\!\!\!\!\!\!& &\!\!\!\!\!\!\horizontal{\circ}\cr
     \horizontal{\circ}\!\!\!\!\!\!& &\!\!\!\!\!\!\horizontal{\symbol}\cr
     \horizontal{\bullet}\!\!\!\!\!\!& &\!\!\!\!\!\!\horizontal{\bullet}\cr
 p\sim -\bar\varrho z &\vertikal{\bullet}\!\!\!\vertikal{\circ}\!\!\!\vertikal{\symbol}\!\!\!\vertikal{\bullet}&p\sim -z\cr
 \end{array}\qquad
\begin{array}{ccc}
p\sim -\varrho z &\vertikal{\circ}\!\!\!\vertikal{\bullet}\!\!\!\vertikal{\symbol}\!\!\!\vertikal{\circ}\!\!\!\vertikal{\bullet}&p\sim -\bar\varrho z\cr
 \horizontal{\circ}\!\!\!\!\!\!& &\!\!\!\!\!\!\horizontal{\bullet}\cr
    \horizontal{\bullet}\!\!\!\!\!\!& &\!\!\!\!\!\!\horizontal{\circ}\cr
     \horizontal{\symbol}\!\!\!\!\!\!& &\!\!\!\!\!\!\horizontal{\symbol}\cr
     \horizontal{\circ}\!\!\!\!\!\!& &\!\!\!\!\!\!\horizontal{\bullet}\cr
 \horizontal{\bullet}\!\!\!\!\!\!& &\!\!\!\!\!\!\horizontal{\circ}\cr
 p\sim -\bar\varrho z &\vertikal{\bullet}\!\!\!\vertikal{\circ}\!\!\!\vertikal{\symbol}\!\!\!\vertikal{\bullet}\!\!\!\vertikal{\circ}&p\sim-\varrho z\cr
 \end{array}
 $$
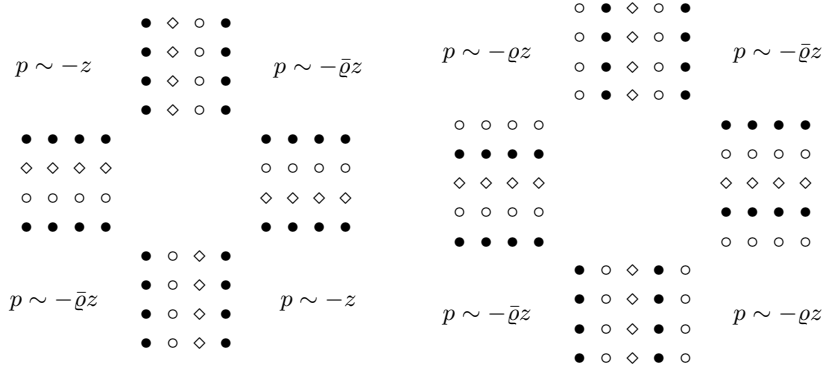
\begin{center}{\sc Figure~4.} Case (A): $\n^\varrho=8$, $\n^{\bar\varrho}=\n^1=4$, and case (B): $\n^\varrho=\n^{\bar\varrho}=8$, $\n^1=4$.
\end{center}}

\section{Application to fourth Painlev\'e transcendents}

\subsection{Fourth Painlev\'e transcendents with positive deficiency $\delta(0,w)$}It is well-known and easy to prove that any transcendental
solution to equation (\ref{PIVmod}) assumes every value $c\ne 0$ equally often in the sense that
$m(r,1/(w-c))=O(\log r).$ This is also true if $c=0$, but $\alpha-\beta+1\ne 0$. In case of $\alpha-\beta+1=0$ there are, however, solutions with Picard value zero; they satisfy $T(r,w)=O(r^2)$. Ignoring this case, it is only known that $m(r,1/w)\le\frac 12T(r,w)+O(\log r)$ holds (equivalently, $m(r,1/y)\le\frac 12T(r,y)+O(\log r)$, see \cite{GLS,NStWertvert}). We will prove the following surprising and never expected result:

\begin{theo}\label{DEFEKT}Let $w=p+q-z$ be any transcendental solution to (\ref{PIVmod}) with positive Nevanlinna deficiency
$\delta(0,w)=1-\limsup\limits_{r\to\infty}\frac{N(r,1/w)}{T(r,w)}.$
Then $(p,q)$ belongs to $\P_{\alpha_k,\alpha_k+1}$ with $\alpha_k=\frac12(-1\pm (2k+1)\sqrt 3 i)$, $k\in\N_0$, and $\delta(0,w)=\frac1{2k+1}$.
For $k=0$, $w$ satisfies the Riccati equation
$w'=\frac{\mp i}{\sqrt 3} (3zw+w^2)$, and zero is a Picard value.\end{theo}

\proof Since $w=p+q-z$ vanishes at least twice at poles of $p$ with residue $1$, Theorem~\ref{EQUALN} yields
$$\ts T(r,w)=\frac23 N(r,\Lambda)+O(r^2)=2N(r,\Lambda^1)+O(r^2)\le N(r,1/w)+O(r^2).$$
This proves $\delta(0,w)=0$ provided $\ds\limsup\limits_{r\to\infty}{T(r,w)}/r^{2}=\infty.$
Otherwise $T(r,w)=O(r^2)$ holds, and $(p,q)$ belongs to the class $\P_{\alpha,\alpha+1}$ for some $\alpha$. From (\ref{Hsystem}) and $\alpha-\beta+1=0$ it follows that
${w'}/w=p-q$, hence  $w$ has zeros of order two at the poles of $p$ with residue $1$, and no others. In any case we have
\be{DELTA0}1-\delta(0,w) =\frac{2\n^1(p)}{\n^\varrho(p)+\n^{\bar\varrho}(p)},\ee
and $\delta(0,w)>0$ implies $\n^1(p)<\frac12(\n^\varrho(p)+\n^{\bar\varrho}(p))$.
It suffices to discuss the signatures (\ref{n1klein}), thus to assume $\n^\varrho\ge\n^{\bar\varrho}$ and $\n^\varrho>\n^1$.
In both cases, the proof of Theorem~\ref{TK} shows that the sequence
{\small$$(p,q){\buildrel{\sf M}_\varrho\over\longrightarrow}(p_1,q_1){\buildrel{\sf B}_1\over\longrightarrow}(p_2,q_2){\buildrel{\sf M}_\varrho\over\longrightarrow}(p_3,q_3){\buildrel{\sf B}_1\over\longrightarrow}(p_4,q_4){\buildrel{\sf M}_\varrho\over\longrightarrow}(p_5,q_5)\cdots$$}
terminates at $(p_n,q_n)\in\P_{\alpha_n,\alpha_n+1}(1)$. We obtain
{\small$$\begin{array}{ccccc}
\schema{\bar\varrho}1{\bar\varrho}1\quad{\buildrel{\sf M}_\varrho\over\longrightarrow}&\schema\varrho{\bar\varrho}\varrho{\bar\varrho}\quad{\buildrel{\sf B}_1\over\longrightarrow}&
\schema{\bar\varrho}\varrho{\bar\varrho}\varrho\quad{\buildrel{\sf M}_\varrho\over\longrightarrow}&\schema{\varrho}1{\varrho}1\quad{\buildrel{\sf B}_1\over\longrightarrow}&\schema{\bar\varrho}1{\bar\varrho}1\quad\cdots\cr
\n^\varrho>\n^{\bar\varrho}=\n^1&\n^1>\n^\varrho=\n^{\bar\varrho}&\n^\varrho=\n^{\bar\varrho}>\n^1&\n^1=\n^\varrho>\n^{\bar\varrho}&\n^\varrho>\n^{\bar\varrho}=\n^1\cr\end{array}
\leqno{\rm (A)}$$}%
and
{\small$$\begin{array}{ccccc}
\schema{\bar\varrho}\varrho{\bar\varrho}\varrho\quad{\buildrel{\sf M}_\varrho\over\longrightarrow}&\schema{\varrho}1{\varrho}1\quad{\buildrel{\sf B}_1\over\longrightarrow}&\schema{\bar\varrho}1{\bar\varrho}1\quad{\buildrel{\sf M}_\varrho\over\longrightarrow}&\schema\varrho{\bar\varrho}\varrho{\bar\varrho}\quad{\buildrel{\sf B}_1\over\longrightarrow}&\schema{\bar\varrho}\varrho{\bar\varrho}\varrho\quad\cdots\cr
\n^\varrho=\n^{\bar\varrho}>\n^1&\n^1=\n^\varrho>\n^{\bar\varrho}&\n^\varrho>\n^{\bar\varrho}=\n^1&\n^1>\n^\varrho=\n^{\bar\varrho}&\n^\varrho=\n^{\bar\varrho}>\n^1\cr\end{array}\leqno{\rm (B)}$$}
respectively, and $\alpha_{4k}=\alpha+k\varrho-k$, $\beta_{4k}=\beta-k\bar\varrho+k$, hence $\alpha_{4k+1}=\bar\varrho\alpha+k-k\bar\varrho$, $\beta_{4k+1}=\varrho\beta-k+k\varrho$, $\alpha_{4k+3}=\beta-(k+1)\bar\varrho+k$, $\beta_{4k+3}=\alpha+(k+1)\varrho-k$, and
$$\begin{array}{rcl}
\alpha_{4k+1}-\beta_{4k+1}+1&=&\bar\varrho\alpha-\varrho\beta+1+3k\cr
\alpha_{4k+3}-\beta_{4k+3}+1&=&-(\alpha-\beta+1)+3k+3.
\end{array}$$
in both cases. The procedure stops when $\n^\varrho=\n^{\bar\varrho}=0$, this implying $n=4k+3$ in case (B). Then, however, $\alpha-\beta+1=0$ is impossible.
In case (A) the procedure stops when $n=4k+1$. Then $\alpha-\beta+1=\alpha_{4k+1}-\beta_{4k+1}+1=0$
has the unique solution $\alpha=-\frac12(1+(2k+1)\sqrt 3 i)$, $\beta=\alpha+1$.
In case of $k=0$ $(\alpha=\bar\varrho)$, $w$ satisfy the Riccati equation with the upper sign, and has Picard value zero, and in any case the procedure  may be reversed. It suffices to describe the first non-trivial step ($k=1$): we start with any generic $(p_0,q_0)\in\P_{\alpha_0,\alpha_0+1}(1)$ with $\alpha_0=\frac12(-1+3\sqrt 3 i)$ and $(\n^1,\n^\varrho,\n^{\bar\varrho})=(4,0,0)$,
and consider the sequence
{\small$$\begin{array}{llllll}
\schema\varrho{\bar\varrho}\varrho{\bar\varrho}\hspace{2.5mm}{\buildrel{\sf M}_{\bar\varrho}\over\longrightarrow}&\schema{\bar\varrho}1{\bar\varrho}1\hspace{2.5mm}{\buildrel{\sf B}_1\over\longrightarrow}
&\schema{\varrho}1{\varrho}1\hspace{2.5mm}{\buildrel{\sf M}_{\bar\varrho}\over\longrightarrow}&\schema{\bar\varrho}\varrho{\bar\varrho}\varrho\hspace{2.5mm}{\buildrel{\sf B}_1\over\longrightarrow}&\schema\varrho{\bar\varrho}\varrho{\bar\varrho}\hspace{2.5mm}{\buildrel{\sf M}_{\bar\varrho}\over\longrightarrow}&\schema{\bar\varrho}1{\bar\varrho}1\cr
(4,0,0)&(0,4,0)&(4,4,0)&(0,4,4)&(8,4,4)&(4,8,4)\cr
\end{array}$$}%
Then $w=p+q-z$ with $(p,q)={\sf M}_{\bar\varrho}{\sf B}_1{\sf M}_{\bar\varrho}{\sf B}_1{\sf M}_{\bar\varrho}(p_0,q_0)\in\P_{-\frac12(1+3\sqrt 3i),\frac12(1-3\sqrt 3i)}$ has deficiency
$$\delta(0,w)=\frac{\n^\varrho+\n^{\bar\varrho}-2\n^1}{\n^\varrho+\n^{\bar\varrho}}=\frac{8+4-2\cdot 4}{8+4}=\frac 13.$$
In the next step ($k=2$ and $\alpha_0=\frac12(-1+5\sqrt 3i)$) we would land at $(p,q)\in\P_{-\frac12(1+5\sqrt 3i),\frac12(1-5\sqrt 3i)}$ with $\n^\varrho=12$, $\n^{\bar\varrho}=\n^1=8$ and $\delta(0,w)=\frac15$ etc. This is also true in the degenerate cases. \Ende

\medskip{\sc Remark.} It is not hard to prove that the exceptional solutions satisfy some implicit  first order algebraic differential equation $P(z,w,w')=0$ of degree $2k+1$ with respect to $w'$.
The corresponding parameters in (\ref{PIV}) are $\hat\alpha=2k+1$ and $\hat\beta=0$.

\subsection{An open problem}The transcendental solutions to (\ref{Hsystem}) as well as to (\ref{PIVmod}) and (\ref{PIV}) have Nevanlinna characteristic
bounded below by $C_1r^2$ and above by $C_2r^4$. It is still open whether or not  {\it either $T(r,p)\asymp r^2$ or else} $T(r,p)\asymp r^4$ holds.

\subsection{Acknowledgement.} I would like to thank the referee for his helpful comments.

\bigskip\noindent{\footnotesize Norbert Steinmetz, Fakult\"at f\"ur Mathematik, Technische Universit\"at Dortmund\\
{\sf stein@math.tu-dortmund.de,
http://www.mathematik.tu-dortmund.de/steinmetz/}}

\end{document}